\documentclass[11pt,leqno]{article}

\usepackage{amsmath,amsfonts,amscd,amssymb,theorem}

\long\def\comment#1\endcomment{}

\comment
\endcomment

\makeatletter
\begingroup
\gdef\th@dotted{\normalfont\itshape
  \def\@begintheorem##1##2{%
        \item[\hskip\labelsep \theorem@headerfont ##1\ ##2.]}%
\def\@opargbegintheorem##1##2##3{%
   \item[\hskip\labelsep \theorem@headerfont ##1\ ##2\ (##3).]}}
\endgroup
\makeatother

\theoremstyle{dotted}

\newtheorem{theorem}{Theorem}[section]
\newtheorem{lemma}[theorem]{Lemma}

\newtheorem{prop}[theorem]{Proposition}
\newtheorem{corr}[theorem]{Corollary}

\makeatletter
\begingroup
\gdef\th@upshape{\normalfont
  \def\@begintheorem##1##2{%
        \item[\hskip\labelsep \theorem@headerfont ##1\ ##2.]}%
\def\@opargbegintheorem##1##2##3{%
   \item[\hskip\labelsep \theorem@headerfont ##1\ ##2\ (##3).]}}
\endgroup
\makeatother

\theoremstyle{upshape}

\newtheorem{defn}[theorem]{Definition}
\newtheorem{remark}[theorem]{Remark}
\newtheorem{exa}[theorem]{Example}

\makeatletter
\renewcommand{\subsection}{\@startsection{subsection}{2}{0pt}{-3ex
plus -1ex minus -0.2ex}{-2mm plus -0pt minus
-2pt}{\normalfont\bfseries}} 
\renewcommand{\subsubsection}{\@startsection{subsubsection}{3}{0pt}{-3ex
plus -1ex minus -0.2ex}{-2mm plus -0pt minus
-2pt}{\normalfont\bfseries}} 
\makeatother

\makeatletter
\@addtoreset{equation}{section}
\makeatother

\newcommand{\cntrct}                
{\hspace{2pt}\raisebox{1pt}{\text{$\lrcorner$}}\hspace{2pt}}

\newcommand{\proof}[1][Proof.]{\smallskip\noindent{\em #1}}
\def\endproof{\hfill\ensuremath{\square}\par\medskip}

\renewcommand{\labelenumi}{{\normalfont(\roman{enumi})}}

\def\eqref#1{\thetag{\ref{#1}}}

\let\latexref=\ref
\def\ref#1{{\normalfont{\latexref{#1}}}}

\newcommand{\wt}{\widetilde}

\newcommand{\dg}{\dagger}

\newcommand{\idot}{{\:\raisebox{1pt}{\text{\circle*{1.5}}}}}
\newcommand{\hdot}{{\:\raisebox{3pt}{\text{\circle*{1.5}}}}}
\newcommand{\Z}{{\mathbb Z}}

\newcommand{\eps}{\varepsilon}
\renewcommand{\phi}{\varphi}

\newcommand{\vH}{\check{H}}

\def\dlim_#1{{\displaystyle\lim_{#1}}}

\newcommand{\Hom}{\operatorname{Hom}}

\newcommand{\Coker}{\operatorname{Coker}}
\newcommand{\Ker}{\operatorname{Ker}}
\renewcommand{\Im}{\operatorname{Im}}

\newcommand{\Fun}{\operatorname{Fun}}

\newcommand{\id}{\operatorname{\sf id}}
\newcommand{\Id}{\operatorname{\sf Id}}
\newcommand{\gr}{\operatorname{\sf gr}}
\newcommand{\tr}{\operatorname{\sf tr}}

\newcommand{\A}{\mathcal{A}}
\newcommand{\D}{\mathcal{D}}
\newcommand{\C}{\mathcal{C}}

\newcommand{\hush}{\natural}

\newcommand{\Aut}{{\operatorname{Aut}}}

\newcommand{\Add}{\operatorname{\sf Add}}

\newcommand{\amod}{{\text{\rm -mod}}}

\newcommand{\ppt}{{\sf pt}}

\newcommand{\lotimes}{\overset{\sf\scriptscriptstyle L}{\otimes}}

\newcommand{\Spec}{\operatorname{Spec}}

\newcommand{\cchar}{\operatorname{\sf char}}

\newcommand{\M}{\operatorname{\mathcal{M}}}

\newcommand{\Infl}{\operatorname{\sf Infl}}

\newcommand{\copr}{{\textstyle\coprod}}

\newcommand{\motimes}{\circ}

\newcommand{\Tr}{\operatorname{Tr}}

\newcommand{\E}{\mathcal{E}}

\newcommand{\vpi}{\check{\pi}}

\newcommand{\LR}{\Lambda R}

\newcommand{\calo}{\mathcal{O}}

\newcommand{\bC}{\overline{C}}

\title{Witt vectors as a polynomial functor}

\author{D. Kaledin\thanks{Partially supported by a subsidy granted
    to the HSE by the Government of the Russian Federation for the
    implementation of the Global Competitiveness Program.}}

\date{\em To Sasha Beilinson, on his birthday.}

\begin{document}

\maketitle

\tableofcontents

\section*{Introduction.}

Recall that to any commutative ring $A$, one canonically associates
the ring $W(A)$ of {\em $p$-typical Witt vectors} of $A$ (a reader
not familiar with the subject can find a great overview for example
in \cite[Section 0.1]{I}). Witt vectors are functorial in $A$, and
$W(A)$ is the inverse limit of rings $W_m(A)$ of {\em $m$-truncated
  $p$-typical Witt vectors} numbered by integers $m \geq 1$. We have
$W_1(A) \cong A$, and for any $m$, $W_{m+1}(A)$ is an extension of
$A$ by $W_m(A)$.

If $A$ is annihilated by a prime $p$ and perfect --- that is, the
Frobenius endomorphism $F:A \to A$ is bijective --- then one has
$W_m(A) \cong W(A)/p^m$, and in particular, $A \cong W(A)/p$. If $A$
is not perfect, this is usually not true. However, if $A$ is
sufficiently nice --- for example, if it is the algebra of functions
on a smooth affine algebraic variety --- then $W(A)$ has no
$p$-torsion. Thus roughly speaking, the Witt vectors construction
provides a functorial way to associate a ring of characteristic $0$
to a ring of characteristic $p$.

Historically, this motivated a lot of interest in the
construction. In particular, one of the earliest attempts to
construct a Weil cohomology theory, due to J.-P. Serre \cite{S}, was
to consider $H^\hdot(X,W(\calo_X))$, where $X$ is an algebraic
variety over a finite field $k$ of positive characteristic $p$, and
$W(\calo_X)$ is the sheaf obtained by taking the Witt vectors of its
structure sheaf $\calo_X$.

This attempt did not quite work, and the focus of attention switched
to other cohomology theories discovered by A. Grothendieck: \'etale
cohomology first of all, but also cristalline cohomology introduced
slightly later. Much later, P. Deligne and L. Illusie \cite{I}
discovered what could be thought of as a vindication of Serre's
original approach. They proved that any smooth algebraic variety $X$
over a perfect field $k$ of positive characteristic can be equipped
with a functorial {\em de Rham-Witt complex} $W\Omega_X^\hdot$, an
extension of the usual de Rham complex $\Omega^\hdot_X$. In degree
$0$, one has $W\Omega^0_X \cong W(\calo_X)$, but in higher degrees,
one needs a new construction. The resulting complex computes
cristalline cohomology $H^\hdot_{cris}(X)$ of $X$, in the sense that
for proper $X$, one has a canonical isomorphism $H^\hdot_{cris}(X)
\cong H^\hdot(X,W\Omega^\hdot_X)$, and cristalline cohomology is
known to be a Weil cohomology theory.

Yet another breakthrough in our understanding of Witt vectors
happenned in 1995, and it was due to L. Hesselholt \cite{heWi}. What
he did was to construct Witt vectors $W(A)$ for an arbitrary
associative ring $A$. Hesselholt's $W(A)$ is also the inverse limit
of its truncated version $W_m(A)$, and if $A$ is commutative and
unital, then it coincides with the classical Witt vectors ring. But
if $A$ is not commutative, $W(A)$ is not even a ring --- it is only
an abelian group. We have $W_1(A) = A/[A,A]$, the quotient of the
algebra $A$ by the subgroup spanned by commutators of its elements,
and for any $m$, $W_{m+1}(A)$ is an extension of $A/[A,A]$ by a
quotient of $W_m(A)$.

In the context of non-commutative algebra and non-commutative
algebraic geometry, one common theme of the two constructions is
immediately obvious: Hochschild homology. On one hand, for any
associative ring, $A/[A,A]$ is the $0$-th Hochschild homology group
$HH_0(A)$. On the other hand, for a smooth affine algebraic variety
$X = \Spec A$, the spaces $H^0(X,\Omega^i_X)$ of differential forms
on $X$ are identified with the Hochshchild homology groups $HH_i(A)$
by the famous theorem of Hochschild, Kostant and Rosenberg. Thus one
is lead to expect that a unifying theory would use an associative
unital $k$-algebra $A$ as an input, and product what one could call
``Hochschild-Witt homology groups'' $WHH_\idot(A)$ such that in
degree $0$, $WHH_0(A)$ coincides with Hesselholt's Witt vectors,
while for a commutative $A$ with smooth spectrum $X = \Spec A$, we
would have natural identifications $WHH_i(A) \cong
H^0(X,W\Omega^i_X)$.

But ``expect'' is perhaps a wrong word here, since this relation to
Hoch\-schild homology is already made abundantly clear in
Hesselholt's work. In fact, while the main construction in
\cite{heWi} is purely algebraic, its origins are in algebraic
topology --- specificaly, the theory of Topological Cyclic Homology
and cyclotomic trace of Bokstedt, Hsiang and Madsen (\cite{BHM},
\cite{HM}), itself a development of Topological Hochschild Homology
of Bokstedt \cite{bo}. Among other things, \cite{BHM} associates a
certain spectrum $TR(A,p)$ to any ring spectrum $A$. If $A$ is a
$k$-algebra, $\cchar k = p$, then $TR(A,p)$ is an Eilenberg-Mac Lane
spectrum, thus essentially a chain complex. The Witt vectors group
$W(A)$ constructed explicitly by Hesselholt is the homology group of
this complex in degree $0$. He also proved (\cite{hedRW},
\cite{heCo}) that for a commutative $k$-algebra $A$ with smooth
spectrum $X$, all the homotopy groups of the spectrum $TR(A,p)$
coincide with the de Rham-Witt forms $H^0(X,W\Omega^\hdot_X)$. Thus
the homotopy groups of $TR(A,p)$ are already perfectly good
candidates for hypothetical Hochschild-Witt homology groups. What we
lack is an algebraic construction of these groups in degrees $> 0$.

\medskip

However, while the relevance of Hochschild homology for Witt vectors
has been well-understood by topologists from the very beginning, one
of its features has been somewhat overlooked. Namely, Hochschild
homology is in fact a theory with two variables --- an algebra $A$
and an $A$-bimodule $M$ (that is, a module over the product $A^o
\otimes A$ of $A$ with its opposite algebra $A^o$). To obtain
$HH_\idot(A)$, one takes as $M$ the diagonal bimodule $A$, but the
groups $HH_\idot(A,M)$ are well-defined for any bimodule. Moreover,
Hochschld homology has the following trace-like property: for any
two algebras $A$, $B$, a left module $M$ over $A^o \otimes B$, and a
left module $N$ over $B^o \otimes A$, we have a canonical
isomorphism
$$
HH_\idot(A,M \lotimes_B N) \cong HH_\idot(B,N \lotimes_A M),
$$
subject to some natural compatibility conditions. This was
axiomatized and studied in the recent work of K. Ponto \cite{pon},
with references to even earlier work. It has also been observed and
axiomatized under the name of ``trace theory'' and ``trace functor''
in \cite{Ka.tr}. The essential point is this: if one wants to have a
generalization of Hochschild homology that is a functor of two
variables $A$, $M$ and has trace isomorphisms, then it suffices to
define it for $A=k$ --- there is general machine that automatically
produces the rest. Thus one can trade the first variable for the
second one: instead of constructing $WHH_\idot(A)$ for an arbitrary
$A$, one can construct $WHH_\idot(k,M)$ for an arbitrary $k$-vector
space $M$. This is hopefully simpler. In particular, it is
reasonable to expect that $WHH_i(k,M)=0$ for $i \geq 1$, so that the
problem reduces to constructing a single functor from $k$-vector
spaces to abelian groups.

\medskip

One possible approach to this is to go back to the original paper
\cite{Bl} of S. Bloch that motivated \cite{I} and \cite{hedRW}, and
use algebraic $K$-theory. It seems that this approach indeed works;
in fact, a definition along these lines has been sketched by the
author a couple of years ago and presented in several talks. The
main technical ingredient of the definition is a certain completed
version of $K_1$ of free non-commutative algebras in several
variables. The construction works over an arbitrary commutative
ring, not just over a perfect field, and provides a version of
``universal'', or ``big'' Witt vectors; the $p$-typical Witt vectors
that we discuss here are then extracted by a separate
procedure. Unfortunately, the complete construction turns out to be
longer than one would wish, and the rather elementary nature of the
resulting Witt vectors functor is somewhat obscured by the machinery
needed to define it. This is one of the reasons why despite many
promises, the story has not been yet written down.

The goal of the present paper is to alleviate the situation by
presenting a very simple and direct alternative construction
motivated by recent work of V. Vologodsky. In a nutshell, here is
the basic idea: instead of trying to associate an abelian group to a
$k$-vector space $M$ directly, one should lift $M$ to a free
$W(k)$-module in some way, use it for the construction, and then
prove that the result does not depend on the lifting. The resulting
definition only works over a perfect field $k$ of characteristic $p$
and assumes that we already know the classical Witt vectors ring
$W(k)$. However, it produces directly an inverse system of
$p$-typical Witt vectors functors $W_m$, and it only uses elementary
properties of cyclic groups $\Z/p^n\Z$, $n \geq 0$. The functors
$W_m$ are polynomial, thus the ``polynomial functor'' of the title.

\medskip

The paper is organized as follows. The longish Section~\ref{pr.sec}
contains the necessary preliminaries. Most of the material is
standard, and we include it mostly to set up notation. We also give
a short recollection on the theory of Mackey functors, the only
piece of relatively high technology required in the body of the
paper. In principle, even this could have been avoided. However,
Mackey functors do help --- in particular, they provide for free
some very useful canonical filtrations on our Witt vectors, and
explain the origin of the projection formula relating the Frobenius
and the Verschiebung maps. Subsection~\ref{mack.subs} and
Subsection~\ref{mack.bis.subs} give an overview of Mackey functors,
Subsection~\ref{Q.subs} contains some less standard material
specifically adapted to our needs, and
Subsection~\ref{cycl.mack.subs} shows how things work for cyclic
groups. We should mention that an intimate relation between
$\Z/p^m\Z$-Mackey functors and $m$-truncated $p$-typical Witt
vectors is very well known and documented in the literature.

Having finished with the preliminaries, we give our two main
definitions, those of Witt vectors and of restriction maps between
them. This, together with the proof of correctness, is the subject
of Section~\ref{witt.sec}. Section~\ref{basic.sec} explores the
basic structure of the Witt vectors functors $W_m$ constructed in
Section~\ref{witt.sec}. We also prove that $W_m$ are pseudotensor
functors (this is Subsection~\ref{mult.subs}), and give a slightly
more explicit inductive description of $W_m$ that is closer to
\cite{heWi} and to the original construction of Witt (this is
Subsection~\ref{exp.subs}). Then in the last
Section~\ref{trace.sec}, we recall the definition of a trace functor
from \cite{Ka.tr}, and prove that our Witt vectors $W_m$ are trace
functors in a natural way. We finish the paper with some results on
compatibility between all the structures we have on $W_m$, namely,
the pseudotensor structure of Subsection~\ref{mult.subs}, the trace
functor structure of Subsection~\ref{tr.witt.subs}, and some
remnants of the Mackey functors structure that they carry by
definition.

As mentioned above, after we define Witt vectors functors on
$k$-vector spaces and equip them with trace isomorphisms,
constructing a Hochschild-Witt complex for any $k$-algebra $A$ with
coefficients in an $A$-bimodule $M$ becomes automatic. However, in
the present paper, we do not do this. We feel that exploring the
resulting Hochschild-Witt complex deserves a separate treatment, and
we relegate it to a companion paper \cite{hw}. The same goes for the
comparison results with Hesselholt's Witt vectors and with the de
Rham-Witt complex. Any possible comparison results with $TR(A,p)$
would require much more technology than we presently have, so will
return to it elsewhere. Another very interesting subject for
comparison is a version of non-commutative Witt vectors given
recently by J. Cuntz and Ch. Deninger \cite{CD}. At the moment, we do
not know what is the relation between the two constructions, and we
feel that it deserves further research.

\subsection*{Acknowledgments.}

The paper is the output of a project of some duration, and during
the project, I had an opportunity to benefit from discussions with
many people. However, I must specifically mention at least three of
them. Firstly, it is hard to overestimate how much I owe to
L. Hesselholt who essentially introduced me to the subject and
generously explained many hidden details. Secondly, I am very
grateful to V. Franjou for his enthusiastic support at the early
stages of my work and for a lot of explanations about polynomial
functors. Thirdly, discussions with V. Vologodsky throughout the
project were very useful, and the main idea for the present paper is
essentially due to him.

Last but not least, it is a great pleasure and honor to dedicate the
paper to Sasha Beilinson. Over the years, I have benefited immensely
from talking to him, on all sorts of things both mathematical and
non-mathematical, and his approach to mathematics is one of the main
inspirations for everything I ever did. As they say in Asia, ``sixty
more''!

\section{Preliminaries.}\label{pr.sec}

\subsection{Small categories.}\label{pr.pr.subs}

For any category $\C$, we will denote by $\C^o$ the opposite
category. We will denote by $\ppt$ the point category (one object,
one morphism). For any group $G$, we will denote by $\ppt_G$ the
groupoid with one object with automorphism group $G$. For any integer
$l \geq 1$, we will simplify notation by letting $\ppt_l =
\ppt_{\Z/l\Z}$. We note that for any functor $F:\C_1 \to \C_2$
between categories $\C_1$, $\C_2$, giving a $G$-action on $F$ is
equivalent to extending it to a functor
\begin{equation}\label{F.G}
\wt{F}:\ppt_G \times \C_1 \to \ppt_G \times \C_2
\end{equation}
that commutes with projections to $\ppt_G$. If the functor $F$
admits a left resp.\ right-adjoint functor $F'$, then $G$ acts on
$F'$ by adjunction, and $\wt{F}$ is left resp.\ right-adjoint to
$\wt{F}'$.

For any small category $I$ and any category $\C$, we will denote by
$\Fun(I,\C)$ the category of functors from $I$ to $\C$. For any ring
$A$, we denote by $A\amod$ the category of left $A$-modules, and for
any small category $I$, we will simplify notation by letting
$\Fun(I,A) = \Fun(I,A\amod)$. This is an abelian category; we will
denote by $\D(I,A)$ its derived category. In particular, for any
group $G$, $\Fun(\ppt_G,A)$ is the category $A[G]\amod$ of left
modules over the group algebra $A[G]$. For any functor $\gamma:I \to
I'$ between small categories, we denote by $\gamma^*:\Fun(I',A) \to
\Fun(I,A)$ the pullback functor, and we denote by
$\gamma_!,\gamma_*:\Fun(I,A) \to \Fun(I',A)$ its left and right
adjoint (the left and right Kan extensions along $\gamma$). The
functor $\gamma^*$ is exact, hence descends to derived categories,
and the derived functors $L^\hdot\gamma_!,R^\hdot\gamma_*:\D(I,A)
\to \D(I',A)$ are left and right adjoint to $\gamma^*:\D(I',A) \to
\D(I,A)$.

If the ring $A$ is commutative, then the category $\Fun(I,A)$ is a
symmetric unital tensor category with respect to the pointwise
tensor product. For any functor $\gamma:I' \to I$, the pullback
functor $\gamma^*$ is a tensor functor. Recall that a functor
$F:\C_1 \to \C_2$ between unital symmetric monoidal categories is
{\em pseudotensor} if it is equipped with functorial maps
\begin{equation}\label{eps.mu}
\eps:1 \to F(1), \quad \mu:F(M) \otimes F(N) \to F(M \otimes N),
\quad M,N \in \C_1,
\end{equation}
where $1$ stands for the unit object, and these maps are associative
and unital in the obvious sense. A pseudotensor functor is {\em
  symmetric} if the maps $\mu$ are also commutative. Then a
right-adjoint to a symmetric tensor functor is automatically
symmetric pseudotensor by adjunction. In particular, for any functor
$\gamma:I' \to I$, the functor $\gamma_*:\Fun(I',A) \to \Fun(I,A)$
is symmetric pseudotensor.

We will assume known the notions of a fibration, a cofibration and a
bifibration of small categories originally introduced in
\cite{sga}. We also assume known the following useful base change
lemma: if we are given a cartesian square
$$
\begin{CD}
I'_1 @>{f_1}>> I_1\\
@V{\pi'}VV @VV{\pi}V\\
I' @>{f}>> I
\end{CD}
$$
of small categories, and $\pi$ is a cofibration, then $\pi'$ is a
cofibration, and the base change map $L^\hdot\pi'_! \circ f_1^* \to
f^* \circ L^\hdot\pi_!$ is an isomorphism. Dually, if $\pi$ is a
fibration, then $\pi'$ is a fibration, and $f^* \circ R^\hdot\pi_*
\cong R^\hdot\pi'_* \circ f_1^*$. For a proof, see e.g.\ \cite[Lemma
  1.7]{K}.

\subsection{Trace maps.}

One specific example of a base change situation is a bifibration
$\pi:I' \to I$ whose fiber is equivalent to a groupoid $\ppt_G$. In
this case, for any $E \in \Fun(I',A)$ and any object $i' \in I'$
with image $i=\pi(i') \in I$, the $A$-module $E(i')$ carries a
natural action of the group $G$, and we have base change
isomorphisms
$$
\pi_!E(i) \cong E(i')_G, \qquad \pi_*E(i) \cong E(i')^G,
$$
where in the right-hand side, we have coinvariants and invariants
with respect to $G$. If the group $G$ is finite, then for any
$A[G]$-module $V$, we have a natural trace map
\begin{equation}\label{tr.G}
\tr_G = \sum_{g \in G}g:V_G \to V^G
\end{equation}
whose cokernel is the Tate cohomology group $\vH^0(G,V)$. Taken
together, these maps then define a natural trace map
\begin{equation}\label{tr.gamma}
\tr_\pi:\pi_!E \to \pi_*E.
\end{equation}
This map is functorial in $E$ and compatible with the base
change. If $E = \pi^*E'$ for some $E' \in \Fun(I,A)$, then
\begin{equation}\label{tr.triv}
\tr_\pi = |G|\id:\pi_!\pi^*E' \to \pi_*\pi^*E',
\end{equation}
where $|G|$ is the order of the finite group $G$. In the general
case, we denote by $\tr^\dg_\pi:E \to E$ the composition
\begin{equation}\label{wt.tr}
\begin{CD}
E @>{l}>> \pi^*\pi_!E @>{\pi^*(\tr_\pi)}>> \pi^*\pi_*E @>{r}>> E,
\end{CD}
\end{equation}
where $l$ and $r$ are the adjunction maps, and we denote by
\begin{equation}\label{v.ga}
\vpi_*:\Fun(I',A) \to \Fun(I,A)
\end{equation}
the functor sending $E$ to the cokernel of the map \eqref{tr.gamma}.
For any $i' \in I'$ with $i=\pi(i')$, we have a natural
identification $\vpi_*(E)(i) \cong \vH^0(G,E(i'))$. We note that
even if the categories $I$, $I'$ are not small, the map
$\tr^\dg_\pi:E \to E$ of \eqref{wt.tr} is perfectly well-defined
for every functor $E:I' \to k\amod$.

\begin{lemma}\label{vpi.le}
Assume that the ring $A$ is commutative. Then the functor $\vpi_*$
of \eqref{v.ga} is pseudotensor.
\end{lemma}

\proof{} Since by definition, $\vpi_*$ is a quotient of a
pseudotensor functor $\pi_*$, it suffices to check that the maps
$\mu$ of \eqref{eps.mu} for the functor $\pi_*$ descend to maps
$\vpi_*(M) \otimes_A \vpi_*(N) \to \vpi_*(M \otimes_A N)$, $M,N \in
\Fun(I',A)$. This can be checked pointwise on $I$. Thus we may
assume that $I = \ppt$ is a point, $M$ and $N$ are $A[G]$-modules,
and $\vpi_*$ is the Tate cohomogy functor $\vH^0(G,-)$. Then the
claim is well-known (and easily follows from the obvious equality
$\tr_G(m \otimes n) = \tr_G(m) \otimes n$, $m \in M$, $n \in N^G
\subset N$).
\endproof

Finally, assume that we have a normal subgroup $N \subset G$ with
the quotient $W = G/N$, and a bifibration $\pi:I'' \to I$ with fiber
$\ppt_G$ factors as $\pi = \pi' \circ \pi''$, where $\pi'':I'' \to
I'$ is a bifibration with fiber $\ppt_N$, and $\pi':I' \to I$ is a
bifibration with fiber $W$. Then we have a commutative diagram
\begin{equation}\label{tr.dia}
\begin{CD}
\pi'_! \circ \pi''_! @>{\sim}>> \pi_! @>{\tr_{\pi'}}>> \pi'_*
\circ \pi''_!\\
@V{\pi'_!(\tr_{\pi''})}VV @VV{\tr_\pi}V @VV{\pi'_*(\tr_{\pi''})}V\\
\pi'_! \circ \pi''_* @>{\tr_{\pi'}}>> \pi_* @>{\sim}>> \pi'_* \circ
\pi''_*.
\end{CD}
\end{equation}
In particular, we have
\begin{equation}\label{tr.prod}
\tr_\pi = \pi'_*(\tr_{\pi''}) \circ \tr_{\pi'} = \tr_{\pi'} \circ
\pi'_!(\tr_{\pi''}).
\end{equation}
Moreover, taking cokernels of the vertical maps in \eqref{tr.dia},
we obtain natural maps
$$
\begin{aligned}
l_\pi:&\pi'_! \circ \vpi''_* = \pi'_!(\Coker(\tr_{\pi''})) \cong
\Coker(\pi'_!(\tr_{\pi''})) \to \vpi_*,\\
r_\pi:&\vpi_* \to \Coker(\pi'_*(\tr_{\pi''})) \to
\pi'_*(\Coker(\tr_{\pi''})) = \pi'_* \circ \vpi''_*,
\end{aligned}
$$
such that $r_\pi \circ l_\pi = \tr_{\pi'}$. By adjunction, these
maps induce maps
\begin{equation}\label{lr.pi}
l^\dg_\pi:\vpi''_* \to \pi^{'*} \circ \vpi_*, \qquad
r^\dg_\pi:\pi^{'*} \circ \vpi_* \to \vpi''_*,
\end{equation}
and we have
\begin{equation}\label{fact}
r^\dg_\pi \circ l^\dg_\pi = \tr^\dg_{\pi'}:\vpi''_* \to \vpi_*.
\end{equation}

\subsection{Mackey functors.}\label{mack.subs}

Mackey functors were introduced by A. Dress \cite{dress}, with
current definition due to H. Lindner \cite{lind}. Good expositions
of the theory can be found e.g.\ \cite{may2} and \cite{the}). We
only give a brief overview following \cite[Section 2]{Ka.proma}.

Assume given a group $G$ and a ring $A$. Denote by $\Gamma_G$ the
category of finite $G$-sets --- that is, finite sets equipped with a
left $G$-action. Denote by $Q\Gamma_G$ the category with the same
objects, and with morphisms from $S_1$ to $S_2$ given by isomorphism
classes of diagrams
\begin{equation}\label{domik}
\begin{CD}
S_1 @<{p_1}<< S @>{p_2}>> S_2
\end{CD}
\end{equation}
in $\Gamma_G$, with composition given by pullbacks. Any map $f:S_1
\to S_2$ between $G$-set defines two maps $f_*:S_1 \to S_2$,
$f^*:S_2 \to S_1$ in the category $Q\Gamma_G$, one by letting
$p_1=\id$, $p_2=f$ in \eqref{domik}, and the other by letting
$p_1=f$, $p_2=\id$. Sending $f$ to $f_*$ resp.\ $f^*$ gives
inclusions $\Gamma_G,\Gamma_G^o \subset Q\Gamma_G$. Say that a
functor $E \in \Fun(\Gamma_G^o,A)$ is {\em additive} if for any
$S,S' \in \Gamma_G$ with disjoint union $S \copr S'$, the natural
map
$$
E(S \copr S') \to E(S) \oplus E(S')
$$
is an isomorphism. An {\em $A$-valued $G$-Mackey functor} is a
functor $M$ from $Q\Gamma_G$ to left $A$-modules whose restriction
to $\Gamma_G^o$ is additive. Mackey functors form a full abelian
subcategory $\M(G,A) \subset \Fun(Q\Gamma_G,A)$. The embedding
functor $\M(G,A) \to \Fun(Q\Gamma_G,A)$ admits a left-adjoint
additivization functor $\Add:\Fun(Q\Gamma_G,A) \to \M(G,A)$.

If $G = \{e\}$ is the trivial group consisting only of its unity
element $e$, so that $\Gamma_G=\Gamma$ is the category of finite
sets, then the category $Q\Gamma$ is in fact equivalent to the
category of free finitely generated commutative monoids, and
$\M(G,A)$ is naturally equivalent to $A\amod$. The functor
$\wt{M}:Q\Gamma \to A\amod$ corresponding to an $A$-module $M$ under
the equivalence sends a finite set $S$ to $M[S]$, the sum of copies
of the module $M$ numbered by elements of the set $S$.

In general, any finite $G$-set $S$ decomposes into a disjoint union
of $G$-orbits: we have
\begin{equation}\label{orb.eq}
S = \coprod_{s \in G \backslash S}[G/H_s],
\end{equation}
where $H_s \subset G$ are cofinite subgroups, and $[G/H]$, $H
\subset G$ a cofinite subgroup denotes the quotient $G/H$ considered
as a $G$-set via the action by left shifts. Then by additivity, the
values of a $G$-Mackey functor $M$ at all finite $G$-sets are
completely determined by its values $M([G/H])$, $H \subset G$ a
cofinite subgroup. To determine $M$ itself, one needs to specify
also the maps $f_*:M([G/H_1]) \to M([G/H_2])$, $f^*:M([G/H_2]) \to
M([G/H_1])$ for any $G$-equivariant map $f:[G/H_1] \to [G/H_2]$. For
any composable pair of maps $f$, $g$, we must have $f_* \circ g_* =
(f \circ g)_*$, $g^* \circ f^* = (f \circ g)^*$. In addition, given
two maps $f:[G/H_1] \to [G/H]$, $g:[G/H_2] \to [G/H]$, we must have
\begin{equation}\label{coset}
g^* \circ f_* = \sum_s f_{s*} \circ g_s^*,
\end{equation}
where the sum is over all the terms $[G/H_s]$ in the orbit
decomposition \eqref{orb.eq} of the fibered product $S = [G/H_1]
\times_{[G/H]} [G/H_2]$, and $f_s:[G/H_s] \to [G/H_2]$, $g_s:[G/H_s]
\to [G/H_1]$ are the natural projections. We note that the
existences of the maps $f$, $g$ implies that $H_1,H_2 \subset G$ can
be conjugated to a subgroup in $H \subset G$. If we do the
conjugation, then we have a natural identification
$$
G \backslash S = G \backslash ([G/H_1] \times_{[G/H]} [G/H_2]) \cong
H \backslash ([H/H_1] \times [H/H_2]) \cong H_1 \backslash H / H_2.
$$
Because of this, \eqref{coset} is known as the {\em double coset
  formula}.

\subsection{Fixed points and products.}\label{mack.bis.subs}

For any subgroup $H \subset G$, a finite $G$-set $S$ is also an
$H$-set by restriction, so that we have a natural functor
$\psi^H:\Gamma_G \to \Gamma_H$. This functor preserves pullbacks,
thus induces a functor $Q\psi^H:Q\Gamma_G \to Q\Gamma_H$. The Kan
extension $Q\psi^H_!$ then sends additive functors to additive
functors, so that we obtain a functor
$$
\Psi^H = Q\psi^H_!:\M(G,A) \to \M(H,A).
$$
This is known as the {\em categorical fixed points functor}. Note
that the centralizer $Z_H \subset G$ of the group $H \subset G$ acts
on $\psi^H$, thus on $\Psi^H$, so that $\Psi^H$ can be promoted to a
functor
\begin{equation}\label{wt.psi}
\wt{\Psi}^H:\M(G,A) \to \M(H,A[Z_H]).
\end{equation}
If $H \subset G$ is cofinite, then $\psi^H$ admits a left-adjoint
functor $\gamma^H:\Gamma_H \to \Gamma_G$ that sends a finite $H$-set
$S$ to $(G \times S)/H$, where $H$ acts on $G$ by right shifts. This
functor also preserves pullbacks, and moreover, the induced functor
$Q\gamma^HQ\Gamma_H \to Q\Gamma_G:$ is also adjoint to $Q\psi^H$
(both on the left and on the right). Then by adjunction, we have
\begin{equation}\label{psi.gamma}
\Psi^H \cong Q\gamma^{H*}.
\end{equation}
On the other hand, assume given a normal subgroup $N \subset G$, and
let $W = G/N$ be the quotient group. Then every finite $W$-set is by
restriction a finite $G$-set, so that we obtain a functor $\Gamma_W
\to \Gamma_G$. This has a right-adjoint $\phi^N:\Gamma_G \to
\Gamma_W$ sending a $G$-set $S$ to the fixed points subset
$S^N$. The functor $\phi^N$ preserves pullbacks, and the induced
functor $Q\phi^N_!$ preserves additivity, so that we obtain the {\em
  geometric fixed points functor}
$$
\Phi^N = Q\phi^N_!:\M(G,A) \to \M(W,A).
$$
The functor $\Phi^N$ has a right-adjoint {\em inflation functor}
$$
\Infl^N = \phi^{N*}:\M(W,A) \to \M(G,A).
$$
The inflation functor is fully faithful. Its essential image
consists of $G$-Mackey functors that are {\em supported at $N$} in
the sense that $M([G/H])=0$ unless $H$ contains $N$.

Assume now that the ring $A$ is commutative. The Cartesian product
preserves pullbacks in each variable, this gives a functor
\begin{equation}\label{Q.m}
m:Q\Gamma_G \times Q\Gamma_G \to Q\Gamma_G.
\end{equation}
The Kan extension $m_!$ does {\em not} in any sense preserve
additivity; however, one can still define a product on the category
$\M(G,A)$ by setting
\begin{equation}\label{motimes}
M_1 \motimes M_2 = \Add(m_!(M_1 \boxtimes_A M_2)).
\end{equation}
This is a symmetric unital tensor product. It is right-exact in each
variable. If $G=\{e\}$ is the trivial group, so that $\M(G,A) \cong
A\amod$, then $M_1 \motimes M_2 \cong M_1 \otimes_A M_2$, and the
unit object is the free module $A$. In general, the unit object is
the so-called {\em Burnside Mackey functor} $A$ given by
\begin{equation}\label{A.eq}
\A = \Add(Qp_!A),
\end{equation}
where $p:\Gamma \to \Gamma_G$ is the tautological embedding sending
a finite set to itself with the trivial $G$-action.

Since $\psi^H \circ m \cong m \circ (\psi^H \times \psi^H)$ and
$\phi^H \circ m \cong m \circ (\phi^H \times \phi^H)$, both fixed
points functors $\Psi^H$, $\Phi^H$ are tensor functors with respect
to the product \eqref{motimes}. Since the isomorphism $\psi^H \circ
m \cong m \circ (\psi^H \times \psi^H)$ is $Z_H$-equivariant, the
extended fixed points functor $\wt{\Psi}^H$ of \eqref{wt.psi} is
also a tensor functor.

While the product $M_1 \motimes M_2$ of two $G$-Mackey functors
$M_1$, $M_2$ does not admit a simple description, such a description
is possible for the spaces of maps from $M_1 \motimes M_2 \to M_3$
to a third $G$-Mackey functor $M_3$. Namely, we have the following
useful result.

\begin{lemma}\label{mack.prod.le}
Assume given a commutative ring $A$, a group $G$, and three
$G$-Mackey functors $M_1,M_2,M_3 \in \M(G,A)$. Then giving a map
$$
\mu:M_1 \motimes M_2 \to M_3
$$
is equivalent to giving a map $\mu:M_1([G/H]) \otimes_A M_2([G/H])
\to M_3([G/H])$ for any cofinite subgroup $H \subset G$ such that
for any map $f:[G/H'] \to [G/H]$ of $G$-sets, we have
\begin{equation}\label{mack.mult.eq}
\begin{aligned}
\mu \circ (f_1^* \times f_2^*) &= f_3^* \circ \mu,\\
\mu \circ (f^1_* \times \id) &= f^3_* \circ \mu \circ (\id \times
f_2^*),\\
\mu \circ (\id \times f^2_*) &= f^3_* \circ \mu \circ (f_1^* \times \id),
\end{aligned}
\end{equation}
where $f_1^*$, $f_2^*$, $f_3^*$, resp.\ $f^1_*$, $f^2_*$, $f^3_*$
are the maps $f^*$ resp.\ $f_*$ for the Mackey functors $M_1$,
$M_2$, $M_3$.
\end{lemma}

\proof{} This is a reformulation of \cite[Lemma
  2.6]{Ka.proma}.
\endproof

Informally, the first equation in \eqref{mack.mult.eq} means that
the maps $f^*$ are multiplicative, and the other two equations say
that $f_*$ satisfies a version of the projection formula with
respect to $f^*$. Note that among other things,
Lemma~\ref{mack.prod.le} implies that for any cofinite $H \subset
G$, the evaluation functor $M \mapsto M([G/H])$ has a natural
pseudotensor structure. This can also been seens explicitly as
follows. For $H=G$, we have a functorial identification
$$
M([G/G]) \cong Qp^*(M),
$$
where $p:\Gamma \to \Gamma_G$ is the tautological embedding of
\eqref{A.eq}, so that evaluation at $[G/G]$ is right-adjoint to the
functor $Qp_!$. Since $Qp$ commutes with the product functor
\eqref{Q.m}, $Qp_!$ is tensor, and $Qp^*$ is pseudotensor by
adjunction. For a general cofinite $H$, precompose with the tensor
functor $\Psi^H$.

\subsection{The functor $Q$.}\label{Q.subs}

If the group $G$ is finite, then the trivial subgroup $\{e\} \subset
G$ is cofinite, and the centralizer $\Z_{\{e\}} \subset G$ is the
whole $G$. Denote by
\begin{equation}\label{E.eq}
U = \wt{\Psi}^{\{e\}}:\M(G,A) \to \M(\{e\},A[G]) \cong A[G]\amod
\end{equation}
the corresponding extended fixed points functor \eqref{wt.psi}. By
\eqref{psi.gamma}, $U$ sends a Mackey functor $M \in \M(G,A)$
to its value $M([G/\{e\}])$ at the biggest $G$-orbit $[G/\{e\}]$,
and $G$ acts on $M[G/\{e\}]$ via its action on $G=[G/\{e\}]$ by
right shifts. The functor $U$ has a left and a right-adjoint
$L,R:A[G]\amod \cong \M(\{e\},A[G]) \to \M(G,A)$ given by
\begin{equation}\label{LR.exp}
L(E) = (\psi^*(E))_G, \quad R(E) = (\psi^*(E))^G,
\qquad E \in A[G]\amod,
\end{equation}
where we simplify notation by writing $\psi = \psi^{\{e\}}$, and $G$
acts both on $E$ and on $\psi$. Explicitly, for any $A[G]$-module
$E$ and any subgroup $H \subset G$, we have
\begin{equation}\label{LR.H}
\begin{aligned}
L(E)([G/H]) &\cong (E[G/H])_G \cong E_H,\\
R(E)([G/H]) &\cong (E[G/H])^G \cong E^H.
\end{aligned}
\end{equation}
We note that $E^H$, $E_H$ only depend on $H$ and not on $G$. In
fact, for any subgroup $H \subset G$, the isomorphisms
\eqref{psi.gamma} and \eqref{LR.exp} provide canonical
identifications
\begin{equation}\label{psi.LR}
\Psi^HL(E) \cong L'(E'), \quad \Psi^HR(E) \cong R'(E'), \qquad E \in
A[G]\amod,
\end{equation}
where $L',R':R[H]\amod \to \M(H,A)$ are the functors $L$, $R$ for
the group $H$, and $E'$ is $E$ treated as an $A[H]$-module by
restriction. In particular, taking $H = \{e\}$, we see that $U \circ
L \cong U \circ R \cong \Id$, so that both $L$ and $R$ are fully
faithful. By adjunction, we then have a natural trace map
\begin{equation}\label{tr.ma}
\tr:L \to R.
\end{equation}
For any subgroup $H \subset G$, it is compatible with the identification
\eqref{psi.LR}.

\begin{defn}\label{Q.def}
For any $A[G]$-module $E$, the $G$-Mackey functor $Q(E)$ is the
cokernel of the trace map \eqref{tr.ma}.
\end{defn}

It is clear that $Q(E)$ is functorial in
$E$, so that we obtain a functor
$$
Q:A[G]\amod \to \M(G,A).
$$
Explicitly, for any subgroup $H \subset G$, the map $L(E)([G/H]) \to
R(E)([G/H])$ induced by the trace map \eqref{tr.ma} is the trace map
$\tr_H$ of \eqref{tr.G}, and we have the identification
\begin{equation}\label{Q.tate}
Q(E)([G/H]) \cong \vH^0(G,E[G/H]) \cong \vH^0(H,E),
\end{equation}
where as in Subsection~\ref{pr.pr.subs}, $\vH^\hdot(-,-)$ stands for
Tate cohomology. The maps $f_*$, $f^*$ associated to a
$G$-equivariant map $f:[G/H'] \to [G/H]$ can aslo be described
explicitly, but we will need it only in one case: $H=G$, $H'=N$ is a
normal subgroup with the quotient $W=G/N$, $f:[G/H] \to [G/G]$ is
the only map. Then
\begin{equation}\label{f.lr}
f_* = l^\dg_\pi, \qquad f^* = r^\dg_\pi,
\end{equation}
where $l^\dg_\pi$, $r^\dg_\pi$ are the natural maps \eqref{lr.pi}
for the projection $\pi:\ppt_G \to \ppt$. We also note that for any
subgroup $H \subset G$, \eqref{psi.LR} provides an identification
\begin{equation}\label{psi.Q}
\Psi^HQ(E) \cong Q'(E'), \qquad E \in A[G]\amod,
\end{equation}
where $E'$ is as in \eqref{psi.LR}, and $Q'$ is the functor $Q$ for
the group $H$.

We will also need the following slightly more elaborate description
of the functor $Q$. Let $Q\wt{\psi}$ be the functor \eqref{F.G}
associated to the action of $G$ on $Q\psi$, and consider the full
embedding
$$
A[G]\amod \cong \M(\{e\},A[G]) \subset \Fun(Q\Gamma,A[G]) \cong
\Fun(\ppt_G \times Q\Gamma_G,A).
$$
Then for any $E \in A[G]\amod \subset \Fun(\ppt_G \times
Q\Gamma,A)$, we have
\begin{equation}\label{pi.LR}
L(E) = \pi_!Q\wt{\psi}^*E, \qquad R(E) = \pi_*Q\wt{\psi}^*E,
\end{equation}
where $\pi:\ppt_G \times Q\gamma_G \to Q\Gamma_G$ is the projection,
and the trace map \eqref{tr.ma} coincides with the trace map
$\tr_\pi$ of \eqref{tr.gamma}. Therefore $Q(E) \cong
\vpi_*Q\wt{\psi}^*E$.

\begin{lemma}\label{Q.le}
For any finite group $G$ and commutative ring $A$, the functor
$Q:A[G]\amod \to \M(G,A)$ of Definition~\ref{Q.def} is symmetric
pseudotensor, and for any subgroup $H \subset G$, the functorial
isomorphism \eqref{Q.tate} is compatible with the pseudotensor
structures.
\end{lemma}

\proof{} Since the functor $U$ is tensor, its right-adjoint $R$ is
symmetric pseudotensor by adjunction. Then as in Lemma~\ref{vpi.le},
since $Q$ is a quotient of $R$, the map $\eps$ of \eqref{eps.mu} for
the functor $R$ induces a map $\eps$ for $Q$, and to prove that $Q$
is pseudotensor, it suffices to show that the maps $\mu$ descend to
the corresponding maps for $Q$ --- the associativity, commutivity
and unitality are then automatic. Moreover, compatibility of
\eqref{Q.tate} with the pseudotensor structures is also automatic,
since the second of the isomorphisms \eqref{LR.H} is compatible with
the pseudotensor structures by adjunction.

Since the product \eqref{motimes} is right-exact, $Q(M) \motimes
Q(N)$ is the cokernel of the map
$$
\begin{CD}
(L(M) \motimes R(N)) \oplus (R(M) \motimes L(N)) @>{(\tr \motimes
    \id) \oplus (\id \motimes \tr)}>> R(M) \motimes R(N)
\end{CD}
$$
for any $M,N \in A[G]\amod$, and since $R$ is symmetric, it
suffices to show that there exists a functorial map
\begin{equation}\label{L.bimod}
L(M) \motimes R(N) \to L(M \otimes N)
\end{equation}
that fits into a commutative diagram
$$
\begin{CD}
L(M) \motimes L(N) @>{\tr \motimes \id}>> R(M) \motimes R(N)\\
@VVV @VV{\mu}V\\
L(M \otimes N) @>{\tr}>> R(M \otimes N).
\end{CD}
$$
To see this, use \eqref{pi.LR}. Denote $\gamma=\gamma^{\{e\}}$, and
let $Q\wt{\gamma}$ be the functor \eqref{F.G} associated to the
$G$-action on $Q\gamma$. Since the product of a free $G$-set and an
arbitrary $G$-set is free, we have a commutative diagram
$$
\begin{CD}
Q\Gamma \times \ppt_G \times Q\Gamma_G @>{(Q\wt{\gamma} \times
  \id)\circ(t \times \id)}>>
\ppt_G \times Q\Gamma_G \times Q\Gamma_G @>{\id \times m}>> \ppt_G
\times Q\Gamma_G\\
@V{\id \times Q\wt{\psi}}VV @. @|\\
Q\Gamma \times \ppt_G \times Q\Gamma @>{(\id \times m) \circ (t
  \times \id)}>> \ppt_G
\times Q\Gamma @>{Q\wt{\gamma}}>> \ppt_G \times Q\Gamma_G,
\end{CD}
$$
where $t:Q\Gamma \times \ppt_G \to \ppt_G \times Q\Gamma$ is the
transposition, and since $Q\wt{\psi}^* \cong Q\wt{\gamma}_!$,
this diagram together with \eqref{pi.LR} induces a functorial
projection formula isomorphism
$$
L(E \otimes_A U(M)) \cong L(E) \motimes M, \qquad E \in A[G]\amod, M
\in \M(G,A).
$$
By adjunction, this isomorphism is compatible with the trace maps,
and the inverse isomorphism yields the map \eqref{L.bimod}.
\endproof

\subsection{Cyclic groups.}\label{cycl.mack.subs}

Now fix a prime $p$, and let $G=\Z_p$ be the group of $p$-adic
integers. Then the lattice of cofinite subgroups $H \subset G$ is
very simple --- they are all of the form $p^mG \subset G$, $m \geq
0$. In particular, finite $G$-orbits are numbered by non-negative
integers, and it turns out that the category $\M(G,A)$ admits a
simple explicit description.

Namely, denote by $I$ the groupoid of $G$-orbits and their
isomorphisms, with $[p^m] \in I$ being the orbit
$[G/p^mG]$. Explicitly, $\Aut([p^m]) = \Z/p^m\Z$, so that $I$ is the
disjoint union of groupoids $\ppt_{p^m}$, $m \geq 0$. Let $I_p
\subset I$ be the full subcategory spanned by orbits other than the
point (in other words, $I_p$ is the union of groupoids $\ppt_{p^m}$
with $m \geq 1$). Denote by $i:I_p \to I$ the natural embedding. On
the other hand, for any $m \geq 1$, we have a natural quotient map
$\Z/p^m\Z \to \Z/p^{m-1}\Z$ and the corresponding functor
$\ppt_{p^m} \to \ppt_{p^{m-1}}$. Taking all these functors together,
we obtain a functor
$$
\pi:I_p \to I.
$$
The functor $\pi$ is a bifibration with fiber $\ppt_p$, so that in
particular, we have the trace map $\tr_\pi$ of \eqref{tr.gamma}.

\begin{lemma}\label{cycl.mack.le}
\begin{enumerate}
\item For any ring $A$, the category $\M(\Z_p,A)$ is equivalent to
  the category of triples $\langle E,V,F \rangle$ of a object $E \in
  \Fun(I,A)$ and two maps
$$
\begin{CD}
i^*E @>{V}>> \pi^*E @>{F}>> i^*E
\end{CD}
$$
such that $F \circ V:i^*E \to i^*E$ is equal to the natural map
$\tr^\dg_\pi$ of \eqref{wt.tr}.
\item Assume that the ring $A$ is commutative, and assume given
  $G$-Mackey functors $M_1,M_2,M_3 \in \M(G,A)$ corresponding to
  triples $\langle E_1,V_1,F_1 \rangle$, $\langle E_2,V_2,F_2
  \rangle$, $\langle E_3,V_3,F_3 \rangle$. Then maps $M_1 \motimes
  M_2 \to M_3$ are in a natural one-to-one correspondence with maps
  $\mu:E_1 \otimes_A E_2 \to E_3$ such that $\mu \circ (F_1 \times
  F_2) = F_3 \circ \mu$ and
$$
\mu \circ (V_1 \times \id) = V_3 \circ \mu \circ (\id \times F_2),\quad
\mu \circ (\id \times V_2) = V_3 \circ \mu \circ (F_1 \times \id).
$$
\end{enumerate}
\end{lemma}

\proof{} By definition, we have a natural embedding $e:I \to
\Gamma_G$ that extends further to an embedding $\wt{e}:I \to
Q\Gamma_G$. Fix a map $\rho:[G/p^{m+1}G] \to [G/p^mG]$ for every $m
\geq 0$. Then any morphism $f:[G/p^nG] \to [G/p^mG]$ in $\Gamma_G$
uniquely decomposes as
\begin{equation}\label{rho.deco}
f = \overline{f} \circ \rho^{n-m}, \qquad \overline{f} \in
\Aut([G/p^mG]),
\end{equation}
and all the maps $\rho$ together define a map of functors $\rho:e
\circ i \to e \circ \pi$. In one direction, the equivalence of
\thetag{i} sends a $G$-Mackey functor $M \in \M(G,A)$ to $\wt{e}^*M
\in \Fun(I,A)$, with $V = \rho_*$, $F=\rho^*$; the equality $F \circ
V = \tr^\dg_\pi$ follows from the double coset formula
\eqref{coset}. In the other direction, $\langle E,V,F \rangle$ gives
a $G$-Mackey functor $M$ such that $M([G/p^mG]) = E([p^m])$, $m \geq
0$, and for any map $f:[G/p^nG] \to [G/p^mG]$, we have
$$
f_* = E(\overline{f})^{-1} \circ V^{n-m}, \qquad f^* = F^{n-m} \circ
E(\overline{f}),
$$
where $\overline{f} \in \Aut([p^m])$ is given by
\eqref{rho.deco}. The double coset formula \eqref{coset} then
follows from the equality $F \circ V = \tr^\dg_\pi$. Finally,
\thetag{ii} immediately follows from Lemma~\ref{mack.prod.le}.
\endproof

We note that if for any integer $m \geq 0$ we denote $G_m = G/p^mG
\cong \Z/p^m\Z$, then for any $n \geq m$, we have a fully faithful
inflation functor
\begin{equation}\label{infl.m.n}
\Infl^n_m:\M(G_m,A) \to \M(G_n,A),
\end{equation}
and for any $m \geq 0$, we have the fully faithful inflation functor
\begin{equation}\label{infl.m}
\Infl_m:\M(G_m,A) \to \M(G,A).
\end{equation}
Therefore Lemma~\ref{cycl.mack.le} also describes the category
$\M(G_m,A)$ for any integer $m \geq 0$ --- this is the full
subcategory spanned by triples $\langle E,V,F \rangle$ such that
$E([p^n])=0$ for $n > m$.

Apart from Lemma~\ref{cycl.mack.le}, one can also study $G_m$-Mackey
functors by induction on $m$. Namely, for any $m \geq 1$ and $M \in
\M(G_m,A)$, we have adjunction maps
\begin{equation}\label{lr}
l:L(U(M)) \to M, \qquad r:M \to R(U(M)),
\end{equation}
and we note that since $U(l)$ and $U(r)$ are isomorphisms, both the
kernel $\Ker r$ and the cokernel $\Coker l$ are supported at
$p^{m-1}G_1 \subset G_m$, so that both are effectively
$G_{m-1}$-Mackey functors. We then introduce the following inductive
definition.

\begin{defn}\label{perf.def}
\begin{enumerate}
\item A $G_m$-Mackey functor $M \in \M(G_m,A)$ is {\em perfect} if
  either $m=0$, or the map $l$ of \eqref{lr} is injective, the map
  $r$ is surjective, and both $\Ker r$ and $\Coker l$ are perfect
  $G_{m-1}$-Mackey functors.
\item Assume given a perfect $G_m$-Mackey functor $M$. Then the {\em
  co-standard filtration} $F_\idot$ on $M$ is the increasing filtration
  such that $F_mM=M$ and $F_iM=F_i(\Ker r)$, $0 \leq i \leq m-1$,
  and the {\em standard filtration} $F^\hdot$ on $M$ is the
  decreasing filtration such that $F^mM=\Im l$, and
  $F^iM=q^{-1}(F^i\Coker l)$, $0 \leq i \leq m-1$, where $q:M \to
  \Coker l$ is the quotient map.
\end{enumerate}
\end{defn}

\begin{exa}\label{bad.exa}
Already for $m=1$, Definition~\ref{perf.def}~\thetag{i} is not
vacuous. Indeed, let $k$ be a field of characteristic $p=\cchar k$,
and take $E \in \Fun(I,k)$ with $E([p^m])=0$, $m \neq 1$, and
$E([p])=k$ with the trivial action of $G_1=\Z/p\Z$. Then the trace
map $\tr_{\Z/p\Z}:k \to k$ vanishes, so that $V=F=0$ satisfies the
condition of Lemma~\ref{cycl.mack.le}~\thetag{i}. The resulting
$G_1$-Mackey functor is not perfect.
\end{exa}

As we see from Example~\ref{bad.exa}, not all $G_m$-Mackey functors
are perfect. However, those of interest to us in the rest of the
paper will be, and the standard and co-standart filtrations will
prove useful.

\section{Witt vectors.}\label{witt.sec}

\subsection{Preliminaries on cyclic groups.}

Fix a perfect field $k$ of positive characteristic $p = \cchar
k$. Let $W(k)$ be the ring of $p$-typical Witt vectors of the field
$k$. Since $k$ is perfect, we have $W(k)/p \cong k$, and for any
integer $n \geq 1$, we have the truncated Witt vectors ring $W_n(k)
= W(k)/p^n$. We also have the Frobenius automorphism $F:W(k) \to
W(k)$ lifting the absolute Frobenius automorphism of the field $k$.

As in Subsection~\ref{cycl.mack.subs}, let $G=\Z_p$ be the group of
$p$-adic integers, and for any integer $m \geq 0$, let
$G_m=\Z/p^m\Z=G/p^mG$. Any finite $G$-set $S$ decomposes as
\begin{equation}\label{S.i.copr}
S = \coprod_{i \geq 0}S_{[i]},
\end{equation}
where $S_{[i]} \subset S$ is the union of elements $s \in S$ whose
stabilizer is exactly $p^iG \subset G$. Any finite $G_m$-set $S$ is
canonically a $G$-set via the quotient map $G \to G_m$; its
decomposition \eqref{S.i.copr} only contains terms $X_{[i]}$ with $i
\leq m$.

\begin{lemma}\label{S.i.le}
Assume given integers $m \geq 0$, $n \geq 1$ and a finite $G_m$-set
$S$, and consider the free $W_n(k)$-module $E = W_n(k)[S]$
spanned by $S$, with the $G_m$-action induced from $S$.
\begin{enumerate}
\item We have a natural identification $E^{G_m} = W_n(k)[S/G_m]$.
\item Moreover, assume that $n \geq m$. Then we have a natural
  identification
$$
\vH^0(G_m,E) = \bigoplus_{0 \leq i < m}W_{m-i}(k)[S_{[i]}/G_m],
$$
where $S_{[i]} \subset S$ are the components of the decomposition
\eqref{S.i.copr}.
\end{enumerate}
\end{lemma}

\proof{} \thetag{i} is obvious. For \thetag{ii}, note that
\eqref{S.i.copr} induces a $G_m$-invariant direct sum decomposition
of $E$, so that it suffices to consider the case $S=S_{[i]}$, $0
\leq i \leq m$. In this case, $S$ is actually a $G_i$-set, and the
$G_i$-action on $S$ is free. Therefore the trace map $\tr_{G_m}$ is
equal to $p^{m-i}\tr_{G_i}$ by \eqref{tr.triv}, and the trace map
$\tr_{G_i}$ is a bijection.
\endproof

Note that for any integers $m \geq 0$, $n \geq 1$, the quotient map
$W_{n+1}(k) \to W_n(k)$ induces a natural map
$$
\vH^0(G_m,W_{n+1}(k)[S]) \to \vH^0(G_m,W_n(k)[S]).
$$
Lemma~\ref{S.i.le}~\thetag{ii} then implies that this map is an
isomorphism if $n \geq m$.

Now assume given a free finitely generated $W_n(k)$-module $E$, and
denote by
\begin{equation}\label{V.m}
E_{(m)} = E^{\otimes_{W_n(k)}p^m}
\end{equation}
its $p^m$-th tensor power with the $W_n(k)$-module structure given
by
\begin{equation}\label{F.tw}
a \cdot e = F^m(a)e, \qquad a \in W_n(k), e \in
E^{\otimes_{W_n(k)}p^m}.
\end{equation}
Let the group $G_m$ act on $E_{(m)}$ by permutations. Moreover, let
$E'_{(m)}$ be the same $W_n(k)$-module considered as a
representation of $G_{m+1}$ via the quotient map $G_{m+1} \to G_m$,
and denote
\begin{equation}\label{Q.m.eq}
Q_m(E) = \vH^0(G_m,E_{(m)}), \qquad Q'_m(E) = \vH^0(G_{m+1},E'_{(m)}).
\end{equation}
Then we have
$$
E_{(m)}^{G_m} \cong (E'_{(m)})^{G_{m+1}},
$$
and $\tr_{G_{m+1}}=p\tr_{G_m}$ on $E_{(m)} \cong E'_{(m)}$ by
\eqref{tr.prod} and \eqref{tr.triv}, so that
we obtain a natural map
\begin{equation}\label{r.eq}
r:Q'_m(E) \to Q_m(E).
\end{equation}
Every element $e \in E$ gives an element $e^{\otimes p^m} \in
E'_{(m)}$; this element is $G_{m+1}$-invariant and descends to
elements
\begin{equation}\label{v.m}
e'_{(m)} \in Q'_m(E), \quad e_{(m)} = r(e'_{(m)}) \in Q_m(E).
\end{equation}

\begin{lemma}\label{v.m.le}
Assume that $n \geq m$, and assume given two elements $a,b \in
E$ such that $a = b \mod p$. Then $a'_{(m)} = b'_{(m)}$ and
$a_{(m)} = b_{(m)}$.
\end{lemma}

\proof{} Since $e_{(m)} = r(e'_{(m)})$, it suffices to prove the
first claim.  By functoriality, it suffices to consider the
universal situation: $E=W_n(k)[S]$ is the free module generated by a
set $S$ with two elements, $s_0$ and $s_1$, and we have $a = s_0$,
$b = s_0 + p s_1$. Then elements $\wt{s}$ of the set $S^{p^m}$ are
of the form
\begin{equation}\label{wt.s.f}
\wt{s} = s_{f(1)} \times s_{f(2)} \times \dots \times s_{f(p^m)},
\end{equation}
where $f$ is a function that assigns $0$ or $1$ to any integer $1
\leq j \leq p^m$. For such an element, denote
$$
|\wt{s}| = \sum_{1 \leq j \leq p^m}f(j),
$$
so that $|\wt{s}|$ is the number of integers $j$ with $f(j)=1$. Then
we have
\begin{equation}\label{binom}
b^{\otimes p^m}-a^{\otimes p^m} =
\sum_{\wt{s}}p^{|\wt{s}|}\wt{s} \in E^{(m)} = W_n(k)[S^{p^m}]^{G_m},
\end{equation}
where the sum is over all elements $\wt{s} \in S^{p^m}$ not equal to
$s_0^{p^m}$ --- that is, with $|\wt{s}| \geq 1$. Moreover, if
$\wt{s}$ lies in $S^{p^m}_{(i)}$ for some integer $i$, then
$|\wt{s}|$ must be divisible by $p^{m-i}$, and since $|\wt{s}| \geq
1$, we actually have $|\wt{s}| \geq p^{m-i}$. Since $p^{m-i} \geq
m-i+1$ for any $i \leq m$, Lemma~\ref{S.i.le}~\thetag{ii} implies
that every term in the right-hand side of \eqref{binom} vanishes
after projecting to $\vH^0(G_{m+1},E'_{(m)})$.
\endproof

\subsection{Polynomial Witt vectors.}

Now fix integers $n \geq m \geq 1$, let $E$ be a free
$W_n(k)$-module, and consider the corresponding
$W_n(k)[G_m]$-modules $E_{(m)}$, $E'_{(m-1)}$ of \eqref{V.m}.
Denote
$$
\wt{Q}_m(E) = Q(E_{(m)}), \wt{Q}'_{m-1}(E) = Q(E'_{(m-1)}) \in
\M(G_m,W_n(k)),
$$
where $Q(-)$ is as in Definition~\ref{Q.def}. Note that by
\eqref{Q.tate}, this notation is consistent with \eqref{Q.m.eq} ---
the Mackey functor $\wt{Q}_m(E)$ gives $Q_m(E)$ after evaluation at
the trivial $G_m$-orbit $[G_m/G_m]$, and similarly for
$\wt{Q}'_m(E)$. Note also that since $n \geq m$, every
$W_m(k)$-module is automatically a $W_n(k)$-module by virtue of the
map $W_n(k) \to W_m(k)$, so that we have a full embedding
$W_m(k)\amod \subset W_n(k)\amod$.

\begin{prop}\label{witt.prop}
For any integer $m \geq 1$, there exists functors
$$
\wt{W}_m,\wt{W}'_{m-1}:k\amod \to \M(G_{m-1},W_m(k))
$$
such that for any free $W_n(k)$-module $E$, we have functorial
isomorphisms
$$
\wt{Q}_m(E) \cong \Infl^m_{m-1}\wt{W}_m(E/p), \qquad \wt{Q}'_{m-1}(E) \cong
\Infl^m_{m-1}\wt{W}'_{m-1}(E/p),
$$
where $\Infl$ are the inflation functors \eqref{infl.m.n}.
\end{prop}

\proof{} The proofs of both claims are the same, so let us start
with $\wt{W}_m$. By \eqref{Q.tate}, we know that
$Q(E_{(m)})([G_m/\{e\}])=0$, so that $\wt{Q}_m(E)$ is supported at
$p^{m-1}G_1 \subset G_m$, thus lies in the image of the fully
faithful inlfation functor $\Infl^m_{m-1}$. Denote by $q$ the
functor from free $W_n(k)$-modules to $k$-vector spaces that sends
$E$ to $E/p$. We have to show that $\wt{Q}_m$ canonically factors through
$q$. Since $q$ is essentially surjective, the issue is the
morphisms: we have to show that for two free $W_n(k)$-modules $M$,
$N$, and two maps $a,b:M \to N$ with $q(a)=q(b)$, we have
$\wt{Q}_m(a)=\wt{Q}_m(b)$. Moreover, since the functor $\wt{Q}_m$ obviously
commutes with filtered colimits, it suffices to consider finitely
generated free $W_n(k)$-modules.

Let $E=\Hom_{W_n(k)}(M,N)$. This is also a free finitely generated
$W_n(k)$-module, we have the action map
$$
\alpha:M \otimes_{W_n(k)} E \to N,
$$
and the map $a:M \to N$ decomposes as
$$
\begin{CD}
M = M \otimes_{W_n(k)} W_n(k) @>{\id \otimes \wt{a}}>> N \otimes_{W_n(k)}
E @>{\alpha}>> N,
\end{CD}
$$
where $\wt{a}:W_n(k) \to E$ is the map sending $1$ to $a$.

Now, the $p^m$-th tensor power functor is tensor, and the functor
$Q$ is pseudotensor by Lemma~\ref{Q.le}. Therefore $\wt{Q}_m$ is
pseudotensor, and the map $\wt{Q}_m(a)$ can be decomposed as
$$
\begin{CD}
\wt{Q}_m(M) \cong \wt{Q}_m(M) \motimes \A @>{\id \motimes
  \wt{Q}_m(\wt{a})}>> \wt{Q}_m(M) \motimes \wt{Q}_m(E) @>{\mu}>>\\
@>{\mu}>> \wt{Q}_m(M \otimes_{W_n(k)} E) @>{\wt{Q}_m(\alpha)}>> \wt{Q}_m(N),
\end{CD}
$$
where $\A \in \M(G_m,W_n(k))$ is the Burnside Mackey functor, and
$\mu$ comes from the pseudotensor structure on $\wt{Q}_m$. We also have an
analogous decomposition for $b$, so that in the end, it suffices to
prove that
$$
\wt{Q}_m(\wt{a}) = \wt{Q}_m(\wt{b}):\A \to \wt{Q}_m(E).
$$
But by \eqref{A.eq} and \eqref{Q.tate}, we have $\Hom(\A,\wt{Q}_m(E))
\cong \vH^0(G_m,E_{(m)})$, and in terms of this identification, we
obviously have $\wt{Q}_m(\wt{a})=a_{(m)}$, $\wt{Q}_m(\wt{b})=b_{(m)}$. Then
we are done by the second claim of Lemma~\ref{v.m.le}.

For $\wt{W}'_{m-1}$, the argument is exactly the same, except that
we need to invoke the first claim of Lemma~\ref{v.m.le}.
\endproof

\begin{defn}\label{witt.m.def}
For any $k$-vector space $E$ and integer $m \geq 1$, the {\em
  $m$-truncated extended polynomial Witt vectors Mackey functor}
$\wt{W}_m(E)$ is the $W_m(k)$-valued $G_{m-1}$-Mackey functor
provided by Proposition~\ref{witt.prop}, and the $W_m(k)$-module of
{\em $m$-truncated polynomial Witt vectors}
$$
W_m(E) = \wt{W}_m(E/k)[G_m/G_m]
$$
is its value at the trivial $G_m$-orbit $[G_m/G_m]$.
\end{defn}

We note that by \eqref{Q.tate}, $m$-truncated polynomial Witt
vectors $W_m(E)$ can be also described as follows: we have
\begin{equation}\label{W.m.Q.m}
W_m(E) \cong Q_m(\wt{E}) = \vH^0(G_m,\wt{E}_{(m)}),
\end{equation}
where $\wt{E}$ is any flat $W_m(k)$-module equipped with an
isomorphism $\wt{E}/p \cong E$, and $Q_m$ is the functor
\eqref{Q.m.eq}.

\subsection{Restriction maps.}

The reader will notice immediately that Definition~\ref{witt.m.def}
only uses the functor $\wt{W}_m$ of Proposition~\ref{witt.prop}. The
role of the functor $\wt{W}'_m$ is that it allows one to relate
$m$-truncated Witt vectors for different $m$.

Namely, we have natural maps $r$ of \eqref{r.eq}, and taken
together, they provide a canonical map
\begin{equation}\label{r.W.eq}
r:\wt{W}'_m(E) \to \Infl^m_{m-1}\wt{W}_m(E)
\end{equation}
for any $m \geq 1$ and any $k$-vector space $E$.

On the other hand, if we take $m=0$, so that $G_{m+1}=G_1=\Z/p\Z$,
then we tautologically have $\vH^0(G_1,\wt{E})\cong \wt{E}/p$ for
any $n$ and any $W_n(k)$-module $\wt{E}$ with the trivial
$G_1$-action, so that we have a canonical identification
$\wt{W}'_0(E) \cong E$ (since $G_0=\{e\}$, we have $\M(G_0,W_1(k))
\cong k\amod$, so that $\wt{W}'(E)$ is simply a $k$-vector
space). But then, there also exists a canonical Frobenius-semilinear
identification
$$
\vH^0(\Z/p\Z,E^{\otimes p}) \cong E
$$
constructed e.g.\ in \cite[Lemma 2.3]{K}. Thus for $m=0$, we have an
isomorphism of functors
\begin{equation}\label{c.eq}
\wt{W}_{m+1} \cong \wt{W}'_m.
\end{equation}
It turns out that the same is true for $m \geq 1$, and we will now
prove it.

\begin{defn}
For $n \geq 1$ and any free $W_n(k)$-module $E$, a
$\Z/p\Z$-equivari\-ant Frobenius-semilinear map $c:E \to
E^{\otimes_{W_n(k)} p}$ is {\em admissible} if the induced map
$$
Q(c):\wt{W}'_0(E/p) \to \wt{W}_1(E/p)
$$
is the standard isomorphism \eqref{c.eq}.
\end{defn}

At this point, we do not need to know the precise form of the
identification \eqref{c.eq}. It suffices to say that if
$E=W_n(k)[S]$ is the free module spanned by a set $S$, then the
diagonal map $\delta:S \to S^p$ induces a Frobenius-semilinear map
\begin{equation}\label{c.S}
c_S:E \to E^{\otimes_{W_n(k)} p}, \qquad \sum_s a_s \cdot s \mapsto
\sum_s a_s^p \cdot \delta(s),
\end{equation}
and this map is admissible.

\begin{prop}\label{redu.prop}
Assume given integers $n \geq m \geq 1$, a free $W_n(k)$-module $E$,
and two admissible maps $c_1,c_2:E \to E^{\otimes_{W_n(k)} p}$. Then
the corresponding maps
$$
Q(c_1^{\otimes p^m}),Q(c_2^{\otimes p^m}):\wt{W}'_m(E/p) \to
\wt{W}_{m+1}(E/p)
$$
coincide, and both are isomorphisms.
\end{prop}

\proof{} As in the proof of Proposition~\ref{witt.prop}, it suffices
to consider finitely generated modules. Moreover, it clearly
suffices to consider the case when $E=W_n(k)[S]$ for some set $S$,
and $c_1$ is the standard map $c_S$ of \eqref{c.S}. In this case, at
least $Q(c_1^{\otimes p^m})$ is certainly an isomorphism by
Lemma~\ref{S.i.le}. We have the decomposition \eqref{S.i.copr} of
the product $S^p$, and it induces the decomposition
\begin{equation}\label{p.pm}
S^{p^{m+1}} = \left(S^p\right)^{p^m} = \coprod_f S^p_{[f(1)]} \times
\dots \times S^p_{[f(p^m)]},
\end{equation}
where $f$ is as in \eqref{wt.s.f}. However, if an element $\wt{s}
\in S^{p^{m+1}}$ is fixed by some non-trivial subgroup in $G_{m+1}$,
it must also be fixed by the smallest non-trivial subgroup $\Z/p\Z
\subset G_{m+1}$, the kernel of the quotient map $G_{m+1} \to
G_m$. Therefore in the decomposition \eqref{S.i.copr} for the
$G_{m+1}$-set $S^{p^{m+1}}$, all the terms except for
$S^{p^{m+1}}_{[m+1]}$ lie in the component of \eqref{p.pm}
corresponding to the constant function $f=0$. By Lemma~\ref{S.i.le},
this means that if we decompose
\begin{equation}\label{V.0}
E_{(m+1)} = W_n(k)[S^{p^{m+1}}] \cong E'_{(m)} \oplus E_0,
\end{equation}
where $E_0$ is spanned by all the other components in \eqref{p.pm},
then any element $e \in E_0$ invariant under $G_{m+1}$ vanishes
after projection to the quotient $W_{m+1}(E/p) \cong
\vH^0(G_{m+1},E_{(m+1)})$. The same is also true if we project to
$\vH^0(H,E_{(m+1)})$ for some subgroup $H \subset G_{m+1}$, so that
$Q(E_0)=0$ by \eqref{Q.tate}.

Now decompose
$$
c_2 = b_0 + b_1,
$$
with $b_i$ taking values in $W_n(k)[S^p_{(i)}]$, $i=0,1$, and
consider the binomial decomposition
\begin{equation}\label{binom.eq}
c_2^{\otimes p^m} = \sum_f b_{f(1)} \otimes \dots \otimes
b_{f(p^m)},
\end{equation}
with the same meaning of $f$ as in \eqref{p.pm}. Then all the terms
except for $b_0^{\otimes p^m}$ take values in $E_0$ of \eqref{V.0},
thus vanish after we apply the functor $Q$. Therefore we may assume
right away that $b_1=0$, so that $c_2 = a \circ c_S$ for some
endomorphism $a:E \to E$ of the module $E$. Since $c_2$ is
admissible, we must have $a=\id$ modulo $p$, and then we are done by
Proposition~\ref{witt.prop}.
\endproof

\begin{corr}\label{restr.corr}
For any integer $m \geq 1$, the standard isomorphism of
Proposition~\ref{redu.prop} provides a functorial isomorphism
\eqref{c.eq}, so that the morphisms \eqref{r.W.eq} provide
functorial maps
\begin{equation}\label{R.eq}
R:\wt{W}_{m+1} \to \Infl^m_{m-1}\wt{W}_m, \qquad R:W_{m+1} \to W_m.
\end{equation}
\end{corr}

\proof{} For the first claim, it suffices to prove that the standard
isomorphism of Proposition~\ref{redu.prop} commutes with
$W_\idot(f)$ for any morphism $f:E_1 \to E_2$ of free
$W_n(k)$-modules. If $f$ is a split injection, so that $E_2 \cong
E_1 \oplus E_1'$, then the claim is obvious --- we can choose bases
in $E_1$ and $E_1'$, consider the corresponding base in $E_2$, and
notice that the already the maps \eqref{c.S} commute with $f$. If
$f$ is surjective, the same argument works. A general map is a
composition of a split injection and a surjection (say, via its graph
decomposition).
\endproof

By virtue of Corollary~\ref{restr.corr}, for any $k$-vector space
$E$, we can consider the inverse limit
\begin{equation}\label{ext.W.eq}
\wt{W}(E) = \lim_{\overset{R}{\gets}}\Infl_{m-1}\wt{W}_m(E) \in
\M(G,W(k)),
\end{equation}
where $\Infl_\idot$ are the inflation functors \eqref{infl.m}.
We call it the {\em extended polynomial Witt vectors} $G$-Mackey
functor of the vector space $E$. Evaluating at the trivial orbit
$[G/G]$, we obtain the $W(k)$-module
\begin{equation}\label{W.eq}
W(E) = \wt{W}(E/k)([G/G]) = \lim_{\overset{R}{\gets}}W_m(E).
\end{equation}
We call it the {\em polynomial Witt vectors module} of the vector
space $E$.

\section{Basic properties.}\label{basic.sec}

\subsection{Exact sequences and filtrations.}\label{filt.subs}

Let us now prove some elementary properties of the truncated Witt
vectors functors of Definition~\ref{witt.m.def}. We start with the
following.

\begin{lemma}\label{R.C.le}
For any integer $m \geq 1$, the map $\wt{W}_{m+1} \to \wt{W}_{m+1}$
given by multiplication by $p$ factors as
$$
\begin{CD}
\wt{W}_{m+1} @>{R}>> \Infl^m_{m-1}\wt{W}_m @>{C}>> \wt{W}_{m+1},
\end{CD}
$$
where $R$ is the restriction map \eqref{R.eq}, and $C$ is a certain
functorial map. Moreover, $R$ is surjective, and $C$ is injective.
\end{lemma}

\proof{} Identify $\wt{W}'_m \cong \wt{W}_{m+1}$ as in
Corollary~\ref{restr.corr}, and let $E$ be an arbitrary free
$W_m(k)$-module. Then by \eqref{tr.prod}, \eqref{tr.triv} and
\eqref{pi.LR}, we have a commutative diagram
$$
\begin{CD}
R(E'_{(m)}) @>{\id}>> R(E_{(m)}) @>{p\id}>> R(E'_{(m)})\\
@A{\tr}AA @AA{\tr}A @AA{\tr}A\\
L(E'_{(m)}) @>{p\id}>> L(E_{(m)}) @>{\id}>> L(E'_{(m)}),
\end{CD}
$$
where $\tr$ are the trace maps \eqref{tr.ma}. Taking the cokernels
of these maps, we obtain the desired factorization. The fact that
$R$ is surjective and $C$ is injective then immediately follows from
Lemma~\ref{S.i.le}.
\endproof

Now for any $k$-vector space $E$ and integer $l \geq 1$, denote
\begin{equation}\label{C.l}
C_l(E) = \left(E^{\otimes l}\right)_{\Z/l\Z}, \qquad C^l(E) =
\left(E^{\otimes l}\right)^{\Z/l\Z},
\end{equation}
where the cyclic group $\Z/l\Z$ acts by permutations, and for any $m
\geq 0$, let $C_{(m)}(E) = C_{p^m}(E)$, $C^{(m)}(E)=C^{p^m}(E)$,
with the $k$-vector space structure twisted by the absolute
Frobenius of $k$ as in \eqref{F.tw}. Note that $C_{(m)}(E)$ and
$C^{(m)}(E)$ canonically extend to $G_m$-Mackey functors
$$
\wt{C}_{(m)}(E) = L(E_{(m)}), \qquad \wt{C}^{(m)}(E) = R(E_{(m)}),
$$
in the sense that we have natural identifications
\begin{equation}\label{C.wtC}
\wt{C}_{(m)}(E)([G_m/G_m]) \cong C_{(m)}(E), \quad
\wt{C}^{(m)}(E)([G_m/G_m]) \cong C^{(m)}(E)
\end{equation}
for any $k$-vector space $E$.

\begin{lemma}\label{lr.RC.le}
For any integer $m \geq 1$, the natural maps $R$, $C$ of
Lemma~\ref{R.C.le} fit into functorial short exact sequences
\begin{equation}\label{exa.rc}
\begin{CD}
0 @>>> \wt{C}_{(m)} @>{l}>> \wt{W}_{m+1} @>{R}>>
\Infl^m_{m-1}\wt{W}_m @>>> 0,\\
0 @>>> \Infl^m_{m-1}\wt{W}_m @>{C}>> \wt{W}_{m+1} @>{r}>>
\wt{C}^{(m)} @>>> 0.
\end{CD}
\end{equation}
\end{lemma}

\proof{} As in Lemma~\ref{R.C.le}, identify $\wt{W}_{m+1} \cong
\wt{W}'_m$. Then by definition and \eqref{Q.tate}, we have
$U(\wt{W}'_m(E)) \cong \vH^0(\Z/p\Z,\wt{E}'_{(m)}) \cong E_{(m)}$
for any free $W_m(k)$-module $\wt{E}$ with reduction $E =
\wt{E}/p$. Therefore we can take the adjunctions maps \eqref{lr} as
$l$ and $r$ in \eqref{exa.rc}. Then since $U \circ \Infl^m_{m-1}=0$,
$R \circ l = r \circ C = 0$ by adjunction, and it suffices to prove
that the sequences are exact after evaluation at an arbitrary
$k$-vector space $E$. This immediately follows from
Lemma~\ref{S.i.le} (choose a basis $S$ in $E$, and use the explicit
decompositions of Lemma~\ref{S.i.le}~\thetag{ii} to compute
$\wt{W}_\idot$).
\endproof

\begin{corr}\label{gr.corr}
For any integer $m \geq 1$ and $k$-vector space $E$, the
$G_{m-1}$-Mackey functor $\wt{W}_m(E)$ is perfect in the sense of
Definition~\ref{perf.def}~\thetag{i}, and its assosiated graded
quotients $\gr^i$, $\gr_i$ with respect to the standard and
co-standard filtrations of Definition~\ref{perf.def}~\thetag{ii} are
given by
$$
\gr^i\wt{W}_m(E) \cong \Infl^{m-1}_i\wt{C}_{(i)}(E), \quad
\gr_i\wt{W}_m(E) \cong \Infl^{m-1}_i\wt{C}^{(i)}(E)
$$
for any $ 0 \leq i \leq m-1$.
\end{corr}

\proof{} Clear. \endproof

Let us now evaluate our $G$-Mackey functors at the trivial $G$-orbit
$[G/G]$. Then \eqref{C.wtC} and \eqref{exa.rc} yield functorial
short exact sequences
$$
\begin{CD}
0 @>>> C_{(m)}(E) @>{l}>> W_{m+1}(E) @>{R}>> W_m(E) @>>> 0,\\
0 @>>> W_m(E) @>{C}>> W_{m+1}(E) @>{r}>> C^{(m)}(E) @>>> 0.
\end{CD}
$$
In fact, we can say more. Namely, for any $m \geq 0$ and $k$-vector
space $E$, denote by $\Phi_m(E)$ the image of the trace map
$$
\tr_{G_m}:C_{(m)}(E) \to C^{(m)}(E).
$$
Then for any $m \geq 0$, sending $e \in E_{(m)}$ to $e^{\otimes p}
\in E_{(m+1)}$ gives a functorial $k$-linear map
$$
C:C_{(m)}(E) \to C_{(m+1)}(E).
$$
If the vector space $E$ is finite-dimensional, then we can dualize
this map to obtain a functorial map
$$
R:C^{(m+1)}(E) \to C^{(m)}(E),
$$
and since both sides commute with filtered colimits in $E$, we can
extend this map to arbitrary $k$-vector spaces.

\begin{lemma}\label{Phi.le}
For any $m \geq 0$ and any $k$-vector space $E$, we have functorial
short exact sequences
\begin{equation}\label{exa.cc}
\begin{CD}
0 @>>> C_{(m)}(E) @>{C}>> C_{(m+1)}(E) @>>> \Phi_{m+1}(E) @>>> 0,\\
0 @>>> \Phi_{(m+1)}(E) @>>> C^{(m+1)} @>{R}>> C^{(m)}(E) @>>> 0,
\end{CD}
\end{equation}
and commutative diagrams
\begin{equation}\label{dia.rc}
\begin{CD}
W_{m+1}(E) @>{r}>> C^{(m+1)}(E)\\
@V{R}VV @VV{R}V\\
W_m(E) @>{r}>> C^{(m)}(E)
\end{CD}
\qquad
\begin{CD}
C_{(m)} @>{l}>> W_m(E)\\
@V{C}VV @VV{C}V\\
C_{(m+1)} @>{l}>> W_{m+1}(E).
\end{CD}
\end{equation}
\end{lemma}

\proof{} Note that we have a functorial four-term sequence
$$
\begin{CD}
0 @>>> C_{(m)} @>{C}>> C_{(m+1)} @>{\tr_{G_{m+1}}}>> C^{(m+1)}
@>{R}>> C^{(m)} @>>> 0.
\end{CD}
$$
and Lemma~\ref{S.i.le} immediately shows that the sequence is exact
(indeed, for any vector space $E=k[S]$ with a basis $S$, the
cokernel of the map $\tr_{G_{m+1}}:C_{(m+1)}(E) \to C^{(m+1)}(E)$ is
$\vH^o(G_{m+1},E^{\otimes p^{m+1}})$, this coincides with
$C^{(m+1)}(E)$ by Lemma~\ref{S.i.le}~\thetag{ii}, and dually for the
kernel $\Ker \tr_{G_{m+1}}$). Together with the definition of the
functor $\Phi_\idot$, this yields the exact sequences
\eqref{exa.cc}. As for \eqref{dia.rc}, then since all the maps are
functorial, it suffices to check commutativity after choosing a
basis $S$ in $E$. Then the first claim immediately follows from
Lemma~\ref{S.i.le} and the definition of the restriction map
$R:W_{m+1} \to W_m$, and the second then follows by
Lemma~\ref{R.C.le}.
\endproof

By induction, Lemma~\ref{Phi.le} shows that the functor $C_{(m)}$
has a natural increasing ``co-standard'' $F_\idot$ filtration with
$\gr_iC_{(m)} \cong \Phi_i$, $0 \leq i \leq m-1$, and the functor
$C^{(m)}$ has a natural descreasing ``standard'' filtration
$F^\hdot$ with $\gr^jC^{(m)} \cong \Phi_j$, $0 \leq j \leq m-1$. The
Witt vectors functor $W_m$ has both filtrations, and they are
transversal. Here is a picture of the associated graded quotient
$\gr^\hdot_\idot W_m = \gr^\hdot\gr_\idot W_m = \gr_\idot \gr^\hdot
W_m$:
\begin{equation}\label{table}
\begin{matrix}
&&&&&& \Phi_0\\
&&&&& \Phi_0 & \Phi_1 \\
&&&& \Phi_0 & \Phi_1 & \Phi_2 \\
&&&\hdotsfor{4}\\
&& \Phi_0 & \hdotsfor{3} & \Phi_{m-3}\\
&\Phi_0 & \Phi_1 & \hdotsfor{2} & \Phi_{m-3} & \Phi_{m-2}\\
\Phi_0 & \Phi_1 & \Phi_2 & \hdots & \Phi_{m-3} & \Phi_{m-2} & \Phi_{m-1}
\end{matrix}
\end{equation}
The indices $i$, $j$ correspond to rows and columns of the table
that we number starting from $0$, and we have $\gr^j_i W_m \cong
\Phi_{i+j+1-m}$, or $0$ if $i+j+1 < m$. Multplication by $p$ acts
diagonally and induces an isomorphism $\gr^i_j \cong
\gr^{i-1}_{j+1}$, or vanishes if $j=m-1$. The bottom row is the
subfunctor $C_{(m-1)} \subset W_m$, and the rightmost column is the
quotient functor $C^{(m-1)}$.

When we pass to the inverse limit \eqref{W.eq}, only the standard
filtration $F^\hdot$ survives --- we have $\gr^iW(E) \cong C_{(i)}$,
$i \geq 0$, and $W(E)/F^iW(E) \cong W_i(E)$ for $i \geq 1$. We can
also consider the inverse limit
$$
C^{(\infty)}(E) = \lim_{\overset{R}{\gets}}C^{(m)}(E),
$$
and it also carries the standard filtration. The second exact
sequence of \eqref{exa.rc} then provides a functorial isomorphism
\begin{equation}\label{W.C}
W(E)/p \cong C^{(\infty)}(E),
\end{equation}
and Lemma~\ref{R.C.le} shows that the $W(k)$-module $W(E)$ is
torsion-free.

\subsection{Frobenius and Verschiebung.}

All the cofinite subgroups $p^nG \subset G$, $n \geq 0$ in the group
$G=\Z_p$ are abstractly isomorphic to $G$ itself. For any integer $n
\geq 0$ and coefficient ring $A$, denote by
$$
\Psi^n = \Psi^{p^nG}:\M(G,A) \to \M(G,A)
$$
the categorical fixed points functor with respect to the subgroup
$p^nG \subset G$, and moreover, for any $m \geq n$, denote
$$
\Psi^n=\Psi^{p^nG_{m-n}}:\M(G_m,A) \to \M(G_{m-n},A).
$$
Since by \eqref{psi.gamma} we obviously have $\Psi^n \circ \Infl_m
\cong \Infl_{m-n} \circ \Psi^n$, the notation is consistent. In
fact, it is even consistent in the sense that $\Psi^{n_1} \circ
\Psi^{n_2} \cong \Psi^{n_1 + n_2}$ for any $n_1,n_2 \geq 0$.

\begin{lemma}\label{Psi.le}
For any $k$-vector space $E$ and integers $m \geq n \geq 1$, we have
a functorial isomorphism
$$
\Psi^n\wt{W}_{m+1}(E) \cong \wt{W}_{m-n+1}(E_{(n)}),
$$
where $E_{(n)}$ is the $k$-vector space \eqref{V.m}. These
isomorphisms commute with restriction maps $R$ and induce a
functorial isomorphism
$$
\Psi^n\wt{W}(E) \cong \wt{W}(E_{(n)}).
$$
\end{lemma}

\proof{} For $\wt{W}_m$, the statement immediately follows from its
definition and \eqref{psi.Q}, and compatibility with the restriction
maps is clear from their construction. For $\wt{W}$, pass to the
limit.
\endproof

Lemma~\ref{Psi.le} allows us to evaluate $G$-Mackey functors
$\Infl_m\wt{W}_{m+1}(E)$, $m \geq 0$ and their inverse limit
$\wt{W}(E)$ at non-trivial $G$-orbits $[G/p^nG]$, $n \geq
1$. Namely, for any integer $n \geq 1$, denote
$$
\begin{aligned}
W^n(E) &= \wt{W}(E)([G/p^nG]),\\
W_m^n(E) &= \wt{W}_m(E)([G_{m-1}/p^nG_{m-n-1}]), \qquad m > n,
\end{aligned}
$$
with $W^0(E)=W(E)$ and $W^0_m(E)=W_m(E)$. By \eqref{Q.tate}, we have
a natural isomorphism
\begin{equation}\label{W.m.Q.m.bis}
W^n_m(E) \cong \vH^0(G_{m-n},\wt{E}_{(m)}),
\end{equation}
with the same meaning of $\wt{E}$ as in \eqref{W.m.Q.m}. On the other
hand, Lemma~\ref{Psi.le} provides natural isomorphisms
\begin{equation}\label{cyclo}
W^n(E) \cong W(E_{(n)}), \quad\ W^n_m(E) \cong
W_{m-n}(E_{(n)}), \quad m > n.
\end{equation}
By definition, $W^n(E)$ and $W^n_m(E)$, $m > n$ carry actions of the
group $G_n$. If $m=n+1$, then \eqref{cyclo} reduces to an
isomorphism
\begin{equation}\label{cyclo.bis}
W^n_{n+1}(E) \cong W_1(E_{(n)} \cong E_{(n)}.
\end{equation}
Effectively, $W^n_{n+1}(E)$ is simply $U(\wt{W}_{n+1}(E/k))$, so
that in particular, the isomorphism \eqref{cyclo.bis} is
$G_n$-equivariant. To see what happens for $m \geq m+2$, equip
$W^n_m(E)$ with the standard and co-standard filtrations induced by
those of Definition~\ref{perf.def}~\thetag{ii}. Both filtrations are
indexed by integers $i$, $n \leq i \leq m-1$, and preserved by
$G_n$. Then \eqref{psi.LR} provides isomorphisms
$$
\Psi^n\wt{C}_{(i)}(E) \cong \wt{C}_{(i-n)}(E_{(n)}), \quad
\Psi^n\wt{C}^{(i)}(E) \cong \wt{C}^{(i-n)}(E_{(n)}), \qquad i \geq
n,
$$
and Corollary~\ref{gr.corr} then shows that we have
\begin{equation}\label{gr.eq}
\gr^iW^n_m(E) \cong C_{(i-n)}(E_{(n)}), \quad
\gr^iW^n_m(E) \cong C_{(i-n)}(E_{(n)})
\end{equation}
for any $i$, $n \leq i < m$. The $G_n$-action in both cases is the
residual action of $G_n = G_i/p^nG_{i-n}$ on $C_{(i-n)}(E_{(n)}) =
(E_{(i)})_{G_{i-n}}$ resp.\ $C_{(i-n)}(E_{(n)}) =
(E_{(i)})_{G_{i-n}}$. In particular, we have
\begin{equation}\label{gr.n}
\begin{aligned}
\gr^iW^n_m(E)_{G_n} &\cong C_{(i)}(E) \cong \gr^iW_m(E),\\
\gr_iW^n_m(E)^{G_n} &\cong C^{(i)}(E) \cong \gr_iW_m(E)
\end{aligned}
\end{equation}
for any $m > i \geq n \geq 0$.

\begin{remark}
As soon as $i > n$, the residual $G_n$-action on
$C_{(i-n)}(E_{(n)})$, $C^{(i-n)}(E_{(n)})$ does {\em not} coincide
with action induced by the natural action on $E_{(n)}$. Explicitly,
the action of the generator $1 \in \Z/p^n\Z = G_n$ is induced by the
map
$$
\sigma_i = \sigma_n^{\otimes p^{i-n}} \circ (\sigma_{i-n} \otimes
\id^{\otimes p^n-1}):E_{(i)} \to E_{(i)},
$$
where $\sigma_n:E_{(n)} \to E_{(n)}$ resp.\ $\sigma_{i-n}:E_{(i-n)}
\to E_{(i-n)}$ are the actions of the generators of the groups $G_n$
resp.\ $G_{i-n}$, and we use the identifications $E_{(i)} \cong
(E_{(n)})_{(i-n)} \cong (E_{(i-n)})_{(n)}$. In fact, $\sigma_i$ is
not even conjugate to $\sigma_n^{\otimes p^{i-n}}$. To see this, it
suffices to consider the case when $E$ is finite-dimensional, choose
a basis, and compute the dimensions of the space of invariants and
coinvariants with respect to the two $G_n$-action, e.g.\ by
Lemma~\ref{S.i.le}.
\end{remark}

Consider now the maps $V$ and $F$ provided by
Lemma~\ref{cycl.mack.le}. Explicitly, these can be described in
terms of the isomorphisms \eqref{W.m.Q.m} and \eqref{W.m.Q.m.bis}
--- by \eqref{f.lr}, we have
\begin{equation}\label{FV.lr}
V = l^\dg_\pi:W^1_m(E) \to W_m(E), \qquad F = r^\dg_\pi:W_m(E) \to
W^1_m(E),
\end{equation}
where $l^\dg_\pi$, $r^\dg_\pi$ are the maps \eqref{lr.pi},
$\pi:\ppt_m \to \ppt$ is the natural projection, and we factorize it
as
$$
\begin{CD}
\ppt_m @>{\pi''}>> \ppt_1 @>{\pi'}>> \ppt.
\end{CD}
$$
Iterating the maps $V$ and $F$, we obtain natural $G_n$-invariant
maps
$$
V^n:W^n(E) \to W(E), \qquad F^n:W(E) \to W^n(E),
$$
and similarly for $W^n_m(E)$, $m > n$. Being $G_n$-invariant,
these maps factor through maps
\begin{equation}\label{bar.FV}
\overline{V}^n:W^n(E)_{G_n} \to W(E), \qquad \overline{F}^n:W(E)
\to W^n(E)^{G_n},
\end{equation}
and again similarly for $W^n_m(E)$, $m > n$.

\begin{lemma}\label{V.R.le}
For any $k$-vector space $E$ and positive integers $m > n \geq
  1$, the maps \eqref{bar.FV} fit into short exact sequences
$$
\begin{CD}
0 @>>> W^n_m(E)_{G_n} @>{\overline{V}^n}>> W_m(E) @>{R^{m-n}}>>
W_{m-n}(E) @>>> 0,\\
0 @>>> W_{m-n}(E) @>{C^{m-n}}>> W_m(E) @>{\overline{F}^n}>>
W^n_m(E)^{G_n} @>>> 0,
\end{CD}
$$
where $R$ and $C$ are the maps of Lemma~\ref{R.C.le}.
\end{lemma}

\proof{} Since $R$ and $C$ are maps of Mackey functors, they commute
with $V$ and $F$, and this immediately implies that $R^{n-m} \circ
V^n = V^n \circ R^{n-m} = 0$ on $W^n_m$ and $F^n \circ
C^{m-n}=C^{m-n} \circ F^n = 0$ on $W_{m-n}$. To prove that the
sequences are exact, use induction on $m-n$. If $m=n$, we let
$W^m_m(E)=0$, and the statement is trivially true. For $m > n$, the
lowest non-trivial term in the standard filtration on $W^n_m(E)$ is
$F^mW^n_m(E) = \gr^mW^n_m(E)$, and by \eqref{gr.eq}, we have
$W^n_m(E)/\gr^mW^n_m(E) \cong W^n_{m-1}(E)$. Then since taking
$G_n$-coinvariants is a right-exact functor, the sequence
$$
\begin{CD}
\gr^mW^n_m(E)_{G_n} @>>> W^n_m(E)_{G_n} @>>> W^n_{m-1}(E)_{G_n} @>>> 0
\end{CD}
$$
is exact on the right. However, by \eqref{gr.n}, the map
$\gr^mW^n_m(E)_{G_n} \to W_m(E)$ induced by $\overline{V}^n$ is
injective, so that the sequence is also exact on the left. Then
$W^n_m(E)/\gr^mW^n_m(E)_{G_n} \cong W^n_{m-1}(E)_{G_n}$ and
$W_m(E)/\gr^mW^n_m(E)_{G_n} \cong W_m/\gr^mW_m(E) \cong W_{m-1}(E)$,
so that proving that the first one of our two sequences is exact for
$W_m$ is equivalent to proving it for $W_{m-1}$. By a dual argument,
exactly the same holds for the second sequence; this gives the
induction step.
\endproof

\begin{exa}\label{k.exa}
If $E=k$ is one-dimensional, then it is immediately obvious from
\eqref{W.m.Q.m} that $W_m(k)$ is the free $W_m(k)$-module of rank
$1$, so that our notation is consistent. Then $V:W_{m-1}(k) \to
W_m(k)$, $F:W_m(k) \to W_{m-1}(k)$ are the standard Verschiebung and
Frobenius maps of the Witt vectors ring $W_m(k)$. Lemma~\ref{V.R.le}
simply states that for any $n < m$, $V^n:W_{m-n}(k) \to W_m(k)$ is
injective, with image $p^{m-n}W_m(k) \subset W_m(k)$, and
$F^n:W_m(k) \to W_{m-n}(k)$ is surjective, with kernel $p^nW_m(k)
\subset W_m(k)$.
\end{exa}

\begin{corr}
On $W^n(E)$, the map $\overline{F}^n$ is an isomorphism, while the
map $\overline{V}^n$ fits into a short exact sequence
$$
\begin{CD}
0 @>>> (W^n(E))_{G_n} @>{\overline{V}^n}>> W(E) @>{R}>> W_n(E)
@>>> 0.
\end{CD}
$$
Moreover, for any integer $i$, we have $\vH^{2i+1}(G_n,W^n(E))=0$,
while the group $\vH^{2i}(G_n,W^n(E))$ is canonically isomorphic to
$W_n(E)$.
\end{corr}

\proof{} The first statement follows from Lemma~\ref{V.R.le} by
taking the inverse limit. For the second statement, compute Tate
cohomology $\vH^\hdot(G_n,-)$ by the standard periodic complex, and
note that by Lemma~\ref{cycl.mack.le}, $F^n \circ V^n$ is the trace
map $\tr_{G_n}$.
\endproof

\subsection{Multiplication.}\label{mult.subs}

Since both $k$ and the Witt vectors rings $W_m(k)$, $m \geq 1$ are
commutative, the categories $k\amod$ and $\M(G_{m-1},W_m(k))$, $m
\geq 1$ are symmetric tensor categories. We then have the following
result.

\begin{prop}\label{mult.prop}
For every integer $m \geq 1$, the extended Witt vectors functor
$\wt{W}_m$ of Definition~\ref{witt.m.def} is symmetric pseudotensor,
and the restriction map \eqref{R.eq} of Corollary~\ref{restr.corr}
is compatible with the pseudotensor structures.
\end{prop}

\proof{} Since $\wt{W}_m$ commutes with filtered colimits, it
suffices to prove both claims after restricting it to
finite-dimensional vector spaces. As in the proof of
Proposition~\ref{witt.prop}, let $q$ be the essentially surjective
reduction functor from the category $W_m(k)\amod^{ff} \subset
W_m(k)\amod$ of free finitely generated $W_m(k)$-modules to the
category $k\amod^f \subset k\amod$ of finite dimensional $k$-vector
spaces. Then for any target category $\C$, the pullback functor
$q^*:\Fun(k\amod^f,\C) \to \Fun(W_m(k)\amod^{ff},\C)$ is fully
faithful, and the same is true for the pullback functor
\begin{multline*}
(q \times q)^*:\Fun(k\amod^f \times k\amod^f,\C) \to\\
\to \Fun(W_m(k)\amod^{ff} \times W_m(k)\amod^{ff},\C).
\end{multline*}
Therefore to construct the maps \eqref{eps.mu} for the functor
$\wt{W}_m$, it suffices to construct them for the functor $\wt{W}_m
\circ q \cong \wt{Q}_m$. But we have $\wt{Q}_m(E) = Q(E_{(m)})$, the
$p^m$-th tensor power functor is symmetric and tensor, and the
functor $Q$ is symmetric pseudotensor by Lemma~\ref{Q.le}.

As for the restriction maps \eqref{R.eq}, what we need to prove is
that they commute with the structure maps $\eps$, $\mu$ of
\eqref{eps.mu}, and this can be checked pointwise, that is, after
evaluating at arbitrary $M,N \in k\amod^f$. Then it suffices to
choose bases in $M$ and $N$, and notice that the standard map $c_S$
of \eqref{c.S} is obviously multiplicative.
\endproof

As an immediate corollary of Proposition~\ref{mult.prop}, we can
pass to the inverse limit with respect to the restriction maps
\eqref{R.eq} and obtain a natural symmeric pseudotensor functors
structure on the extended Witt vectors functor $\wt{W}$ of
\eqref{ext.W.eq}. Another immediate corollary is a symmetric
pseudotensor structure on the polynomial Witt vectors functors
$W_m$, $m \geq 1$ and their inverse limit $W$ obtained by evaluating
at the trivial $G$-orbit $[G/G]$. Since by Lemma~\ref{Q.le}, the
pseudotensor structure on the functor $Q$ is compatible with the
isomorphism \eqref{Q.tate}, the pseudotensor structures on $W_m$, $m
\geq 1$ can also be characterized directly --- these are the only
pseudotensor structures compatible with the isomorphism
\eqref{W.m.Q.m} and the obvious pseudotensor structure on the
functor $Q_m$.

We observe that the unit object in $k\amod$ is the one-dimensional
vector space $k$, and the unit object in $W_m(k)\amod$ is the free
module $W_m(k)$ of rank $1$. Then by Example~\ref{k.exa}, the map
$\eps$ of the pseudotensor structure on the functor $W_m$ is an
isomorphism. For any $k$-vector space $E$, applying the map $\mu$ to
the tautological isomorphism $k \otimes E \cong E$ gives a map
$W_m(k) \otimes_{W_m(k)} W_m(E) \to W_m(E)$, and this map is also
tautologically an isomorphisms. As it happens, these tautologies
already give one thing essentially for free.

\begin{lemma}
For any $k$-vector space $E$ and integer $m \geq 1$, we have
$$
V \circ F = p\id:W_m(E) \to W_m(E).
$$
\end{lemma}

\proof{} Applying Lemma~\ref{cycl.mack.le}~\thetag{ii} to the
tautological isomorphism $k \otimes E \cong E$, we see that it
suffices to prove the claim for $E = k$. But then $F$ and $V$
obviously commute, so that $VF=FV=\tr_{\Z/p\Z}$, and the
$\Z/p\Z$-action on $W_m(k)$ is trivial.
\endproof

Another immediate corollary of the tautological isomoprhism $W_m(k)
\cong W_m(k)$ is the existence of a functorial map
\begin{equation}\label{T.m}
T:E = \Hom_k(k,E) \to W_m(E) = \Hom_{W_m(k)}(W_m(k),W_m(E))
\end{equation}
induced by the functor $W_m$. The map $T$ is not additive --- it is
only a map of sets, a version of the Teichm\"uller representative
map for our polynomial Witt vectors $W_m$. However, it is obviously
multiplicative, in the the sense that
\begin{equation}\label{T.mult}
T(e \otimes e') = \mu(T(e) \otimes T(e')), \qquad e,e' \in E,
\end{equation}
and compatible with restriction maps \eqref{R.eq}, in that $R \circ
T = T$. Passing to the limit, we obtain a functorial map of sets
\begin{equation}\label{T}
T:E \to W(E), \qquad E \in k\amod.
\end{equation}
Spelling out the definition of the functor $W_m$, we see that one
can describe the Teichm\"uller map \eqref{T.m} as follows. For any
$e \in E$, choose a free $W_m(k)$-module $\wt{E}$ so that $E \cong
\wt{E}/p$, lift $e$ to an element $\wt{e} \in E$, and consider the
elements $\wt{e}_{(m)} \in Q_m(\wt{E}) = W_m(E)$, $\wt{e}'_{(m-1)}
\in Q'_{m-1}(\wt{E}) = W'_{m-1}(E)$ of \eqref{v.m}. Then both these
elements only depend on $e$ by Lemma~\ref{v.m.le}, the canonical
isomorphism of Proposition~\ref{redu.prop} sends $\wt{e}'_{(m-1)}$
to $\wt{e}_{(m)}$, and we have
\begin{equation}\label{T.exp}
T(e) = \wt{e}_{(m)} = \wt{e}_{(m-1)}' \in W_m(E) \cong W'_{m-1}(E).
\end{equation}
For a more non-trivial application of Proposition~\ref{mult.prop},
note that for any finite-dimensional $k$-vector space $E$ with the
dual vector space $E^*$, we have the natural pairing map $E \otimes
E^* \to k$.  Together with the pseudotensor structure on $W_m$, it
induces a natural pairing
\begin{equation}\label{pair}
\begin{CD}
W_m(E) \otimes_{W_m(k)} W_m(E^*) @>{\mu}>> W_m(E \otimes E^*)
@>>> W_m(k).
\end{CD}
\end{equation}
We observe that by Lemma~\ref{cycl.mack.le}~\thetag{ii}, this
pairing must be compatible with the maps $F$ and $V$ in the
following sense: we have
\begin{equation}\label{FV.adj}
\langle V(a),b \rangle = V(\langle a,F(b) \rangle) \in W_m(k)
\end{equation}
for any $a \in W_{m-1}(E_{(1)})$, $b \in W_m(E^*)$.

\begin{lemma}
For any $m \geq 1$ and finite-dimensional $k$-vector space $E$,
\eqref{pair} is a perfect pairing, so that $W_m(E)$ and $W_m(E^*)$
are dual modules over the Gorenstein ring $W_m(k)$.
\end{lemma}

\proof{} Let $W_m(E)^* = \Hom_{W_m(k)}(W_m(E),W_m(k))$ be the
$W_m(k)$-module dual to $W_m(E)$. Then the pairing \eqref{pair}
induces a map
$$
\rho:W_m(E^*) \to W_m(E)^*,
$$
and we have to prove that this map is an isomorphism. If $m=1$, we
have $W_1(E^*) \cong E^*$, $W_1(E) \cong E$, so the claim is
clear. In the general case, we note that by Lemma~\ref{S.i.le} and
\eqref{W.m.Q.m}, the source and the target of the map $\rho$ are
$W_m(k)$-modules of the same finite length, so that it suffices to
prove that $\rho$ is injective. Assume by induction that this is
proved for $W_{m-1}$. Then since the Verschiebung map $V:W_{m-1}(k)
\to W_m(k)$ is injective, \eqref{FV.adj} implies that for any $a \in
\Ker \rho$, we have $F(a)=0$. By Lemma~\ref{V.R.le}, this means that
$a$ lies in the image of the map $C^{m-1}$, and by induction, we
have $a=0$.
\endproof

We note by Lemma~\ref{V.R.le} and \eqref{FV.adj}, the perfect
pairing \eqref{pair} interchanges the standard and costandard
filtrations on $W_m$. In terms of the diagram \eqref{table}, one
notes that the functors $\Phi_i$ are obviously self-dual, that is,
we have $\Phi_i(E)^* \cong \Phi_i(E^*)$, and the pairing
\eqref{pair} then corresponds to flipping the table along the main
diagonal.

\medskip

Assume now given an algebra $A$ over $k$. Then since $W_m$ is
symmetric pseudotensor, $W_m(A)$, $m \geq 1$ and $W(A)$ are rings,
and if the algebra $A$ is commutative, these rings are also
commutative.

Unfortunately, at this point our notation stops being consistent:
for any commutative $k$-algebra $A$ different from $k$ itself, our
ring $W_m(A)$ is definitely {\em not} isomorphic to the standard
$m$-truncated $p$-typical Witt vectors ring $W^{st}_m(A)$ of $A$,
and $W(A)$ is not isomorphic to the standard Witt vectors ring
$W^{st}(A)$. Indeed, already the associated graded quotients with
respect to the standard filtration are different: we have $\gr^iW(A)
\cong C_{(i)}(A)$ and $\gr^iW^{st}(A) \cong A$. In fact, our
construction is by its nature relative over the field $k$, so that
$W(E/k)$ would have been perhaps better notation for our polynomial
Witt vectors. In this sense, the standard Witt vectors ring ought to
correspond to $W(A/A)$. Unfortunately, under our present definition
of polynomial Witt vectors this makes no sense, since the
constrution crucially depends on the identification \eqref{c.eq},
and it only holds for perfect $k$. As we have mentioned in the
Introduction, it is possible to give a competely different
definition that does not require this assumption, but this is much
more difficult technically, and we will return to it elsewhere.

\subsection{Explicit descriptions.}\label{exp.subs}

Let us now give some more explicit descriptions of the functors
$W_m$. First of all, for any set $S$ and integer $m \geq 0$, denote
$$
C_{(m)}(S) = S^{p^m}/G_m, \qquad \bC_{(m)} = S^{p^m}_{[m]}/G_m =
C_{(m)}(S) \setminus C_{(m-1)}(S),
$$
where the embedding $C_{(m-1)}(S) \subset C_{(m)}(S)$ is induced by
the diagonal embedding $S^{p^{m-1}} \subset S^{p^m}$. Note that this
is consistent with our earlier notation, in that we have
$C_{(m)}(k[S]) \cong k[C_{(m)}(S)]$. Then the Witt vectors $W_m(E)$
of the $k$-vector space $E=k[S]$ spanned by $S$ are immediately given
by Lemma~\ref{S.i.le} --- we have a natural identification
$$
W_m(E) \cong \bigoplus_{0 \leq i < m}W_{m-i}(k)[\bC_{(i)}(S)].
$$
Now more generally, assume given a $k$-vector space $E$ graded by
$S$, that is, assume that we have
$$
E = \bigoplus_{s \in S}E_s
$$
for some $k$-vector spaces $E_s$, $s \in S$. 

\begin{lemma}\label{gr.S.le}
For any $S$-graded vector space $E$, and for any choice of
splittings $\tau:\C_{(i)}(S) \to S^{p^i}_{[i]}$ of the quotient maps
$S^{p^i}_{[i]} \to \C_{(i)}(S)$, $0 \leq i < m$, we have a natural
identification
\begin{equation}\label{gr.S}
W_m(E) \cong \bigoplus_{0 \leq i < m}\bigoplus_{f \in
  \bC_{(i)}(S)}W_{m-i}(E_{\tau(f)(1)} \otimes \dots \otimes
E_{\tau(f)(p^i)}).
\end{equation}
\end{lemma}

\proof{} For any $i \geq 0$, we have a natural decomposition
$$
E_{(i)} = E^{\otimes p^i} \cong \bigoplus_{f \in S^{p^i}}E_{f(1)} \otimes
\dots \otimes E_{f(p^i)},
$$
so that the map $V^i:W_{m-i}(E_{(i)}) \to W_m(E)$ induces a map
$$
V^i:W_{m-i}(E_{\tau(f)(1)} \otimes \dots \otimes
E_{\tau(f)(p^i)}) \to W_m(E)
$$
for any $f \in \bC_{(i)}(S)$. To prove that the direct sum of these
maps is an isomorphism, it suffices to prove that it becomes an
isomorphism after we pass to the associated graded quotients with
respect to the standard filtration. This immediately follows from
\eqref{gr.n} and an obvious computation for the cyclic power
functors $C_{(i)}(-)$, $0 \leq i < m$.
\endproof

\begin{remark}
The independence of the choice of $\tau$ in Lemma~\ref{gr.S.le} is
not surprising. In fact, as we will see in Section~\ref{trace.sec},
$W_m$ has a natural structure of a trace functor in the sense of
\cite{Ka.tr}. In particular, this provides a natural identification
$$
W_m(E_1 \otimes \dots \otimes E_n) \cong W_m(E_{\sigma(1)} \otimes
\dots \otimes E_{\sigma(n)})
$$
for any $k$-vector spaces $E_1,\dots,E_n$ and any cyclic permutation
$\sigma$ of the set of indices $\{1,\dots,n\}$.
\end{remark}

If a $k$-vector space is graded by integers, then we can lift it to
a $\Z$-graded flat $W_m(k)$-module $\wt{E}$, and then by
\eqref{W.m.Q.m}, $W_m(E)$ inherits a natural grading. However, the
restriction maps $R:W_{m+1}(E) \to W_m(E)$ do not preserve this
grading --- conversely, they multiply it by $p$. To make the graded
consistent with the restriction maps, we need to rescale it by
$p^m$. Thus a natural grading on $W_m(E)$, $m \geq 1$ and on the
limit $W(E)$ is indexed not by integers but by elements $a \in
\Z[p^{-1}]$ in the localization at $p$ of the ring $\Z$. To see this
grading in terms of the decomposition \eqref{gr.S}, denote
$$
|f| = \sum_{j=1}^{p^i}\tau(f)(j)
$$
for any $i \geq 0$ and $f \in \bC_{(i)}(\Z)$. Then the component
$W_m(E)_a \subset W_m(E)$ of degree $a$ is given by
\begin{equation}\label{gr.Z}
W_m(E)_a = \bigoplus_{0 \leq i <
  m}\bigoplus_{|f|=p^ia}W_{m-i}(E_{\tau(f)(1)} \otimes \dots \otimes
E_{\tau(f)(p^i)}).
\end{equation}
Passing to the limit, one obtains an obvious version of \eqref{gr.S}
and \eqref{gr.Z} for the limit Witt vectors functor $W$.

We now observe that even if a $k$-vector space $E$ has no
distinguished basis or grading, we still have the Teichm\"uller
representative maps \eqref{T}, \eqref{T.m}. We recall that these
maps are not additive. Nevertheless, combining them with the powers
$V^m$ of the Verschiebung map $V$, we can still obtain a functorial
surjective map of sets
\begin{equation}\label{wt.T}
\wt{T}:\prod_{m \geq 0}E_{(m)} \to W(E)
\end{equation}
given by
$$
\wt{T}(\langle e_0,e_1,\dots \rangle) = \sum_{m \geq 0}V^m(T(e_m)),
$$
where the sum is convergent with respect to the limit topology on
$W(E)$.

In principle, one can use the map \eqref{wt.T} to obtain an
alternative description of the functor $W(E)$ that is closer to the
classical construction of Witt vectors. We will not do it
completely. However, let us show how one can describe the first
non-trivial extension
$$
\begin{CD}
0 @>>> W^1(E) @>{\overline{V}}>> W(E) @>{R}>> E @>>> 0
\end{CD}
$$
in the standard filtration on $W(E)$. The Teichm\"uller map
\eqref{T} gives a set-theoretic section of the restriction map $R$,
so that we obtain a functorial isomorphism of sets
$$
W(E) \cong W^1(E) \times E,
$$
and in terms of this isomorphism, the abelian group structure on
$W(E)$ is given by
\begin{equation}\label{coc.W}
(e'_0 \times e_0) + (e'_1 \times e_1) = (e_0' + e_1' - c(e_0,e_1))
\times (e_0 \times e_1),
\end{equation}
where $c(-,-)$ is a certain functorial symmetric $2$-cocycle of the
group $E$ with coefficients in $W^1(E) \subset W(E)$. For any $m
\geq 2$, we can project the cocycle $c$ to $W_m(E)$ and obtain a
cocycle
\begin{equation}\label{coc.m}
c:E \times E \to W^1_m(E)
\end{equation}
that gives the group $W_m(E)$.

\begin{lemma}\label{coc.le}
\begin{enumerate}
\item Let $M = \Z[S]$ be the free abelian group generated by the set
  $S=\{s_0,s_1\}$ with two elements, and for any integer $n \geq 1$,
  let $\sigma:S^n \to S^n$ be the cyclic permutation of order
  $n$. Then there exist elements $c_i \in M^{\otimes p^i}$, $i \geq
  1$ such that for any $n \geq 1$, we have
\begin{equation}\label{c.n.eq}
(s_0+s_1)^{\otimes p^n} = s_0^{\otimes p^n} + s_1^{\otimes p^n} +
  \sum_{i=1}^n\sum_{j=0}^{p^i-1}\sigma^j(c_i^{\otimes p^{n-i}}) \in
  M^{\otimes p^n}.
\end{equation}
\item Moreover, for any such set of elements $c_i$, $1 \leq i < m$
  satisfying \eqref{c.n.eq} for $n < m$, the cocycle \eqref{coc.m}
  is given by
$$
c(e_0,e_1) = \wt{T}(0 \times c_1(e_0,e_1) \times c_2(e_0,e_1) \times
c_{m-1}(e_0,e_1), \qquad e_0,e_1 \in E,
$$
where $\wt{T}$ is the map \eqref{wt.T}, and $c_i(e_0,e_1)$ stands
for the image of $c_i$ under the map $M \to E$ sending $s_0$ to
$e_1$ and $s_1$ to $e_1$.
\end{enumerate}
\end{lemma}

\proof{} For \thetag{ii}, it suffices to use the interpretation
\eqref{T.exp} of the Teichm\"uller map in terms of the element
$\wt{e}'_{(m-1)}$, and recall that for any $i$ such that $0 < i <
m$, the map $V^i:Q'_{m-1-i}(\wt{E}_{(i)}) \to Q'_{m-1}(\wt{E})$ is
induced by the trace map
$$
\tr_{G_i}:H^0(G_{m-i},\wt{E}'_{(m-1)}) \to H^0(G_m,\wt{E}'_{(m-1)}).
$$
For \thetag{i}, assume by induction that we already have elements
$c_i$, $1 \leq i < m$ satisfying \eqref{c.n.eq} for $n < m$. Then
the difference
$$
\overline{c}_m=(s_0+s_1)^{\otimes p^m} - s_0^{\otimes p^{m-1}} -
s_1^{\otimes p^{m-1}} -
\sum_{i=1}^{m-1}\sum_{j=0}^{p^i-1}\sigma^j(c_i^{\otimes p^{m-i}})
$$
is invariant under $\sigma$, thus gives an element in
$H^0(G_m,M^{\otimes p^m})$. An element $c_m$
satisfies \eqref{c.n.eq} for $n=m$ if and only if
$$
\overline{c}_m = \sum_{j=0}^{p^m-1}\sigma^j(c_m),
$$
thus to show that it exists, it suffices to check that
$\overline{c}_m$ projects to $0$ in the Tate cohomology group
$\vH^0(G_m,M^{\otimes p^m})$. This cohomology group is certainly
annihilated by $p^m = |G_m|$, thus it does not change if we replace
$M$ with $M \otimes \Z_p$, and then \eqref{W.m.Q.m} provides an
identification
$$
\vH^0(G_m,M^{\otimes p^m}) \cong W_m(E),
$$
where $E = M/p$ is treated as a vector space over the prime field
$\Z/p\Z$. To see that $\overline{c}_m$ projects to $0$ in this
group, it suffices to use the interpretation \eqref{T.exp} of the
map $T$ in terms of the element $\wt{e}_{(m)}$, and apply \thetag{ii}.
\endproof

\begin{remark}
Lemma~\ref{coc.le}~\thetag{i} has {\em a priori} nothing to do with
the functors $W_m$, and might admit an explicit combinatorial proof
(for example, certain explicit polynomials $\delta_i$ are introduced
in \cite{heWi}, and the proof of \cite[Proposition 1.2.3]{heWi}
seems to also prove that they satisfy \eqref{c.n.eq}). Having
obtained the universal polynomials $c_i$ in whatever fashion, one
might try to reconstruct the functors $W_m$ by induction on $m$ via
\eqref{coc.W}. However, one would then also need a $\Z/p\Z$-action
on $W_{m-1}(E^{\otimes p})$ that produces the quotient $W^1_m(E) =
W_{m-1}(E^{\otimes p})_{\Z/p\Z}$, and this requires a structure of a
trace functor that we explore in Section~\ref{trace.sec}. One might
be able to reconstruct this structure explicitly as well, but we did
not pursue this.
\end{remark}

\section{Trace functors.}\label{trace.sec}

\subsection{Definitions and examples.}

A {\em trace functor} from a unital monoidal category $\E$ to a
category $\C$ is a functor $P:\E \to \C$ equipped with functorial
isomorphisms
$$
\tau_{M,N}:P(M \otimes N) \to P(N \otimes M), \qquad M,N \in \E
$$
such that for any object $M \in \E$, we have
$\tau_{1,M}=\tau_{M,1}=\id$, and for any three objects $M,N,L \in
\E$, we have
\begin{equation}\label{M.N.L}
\tau_{L,M,N} \circ \tau_{N,L,M} \circ \tau_{M,N,L} = \id,
\end{equation}
where $1 \in \E$ denotes the unit object, we identify $1 \otimes M
\cong M \cong M \otimes 1$ by means of the unitality isomorphism of
the category $\E$, and $\tau_{A,B,C}$ for any $A,B,C \in \E$ is the
composition of the map $\tau_{A,B \otimes C}$ and the map induced by
the associativity isomorphism $(B \otimes C) \otimes A \cong B
\otimes (C \otimes A)$.

It seems that the notion of a trace functor has been around in some
form at least since 1960-ies. This particular definition is taken
from \cite{Ka.tr}, and it admits a convenient repackaging using
A. Connes' cyclic category $\Lambda$. We recall (see e.g.\ \cite{Lo}
for \cite{FT}) that objects in $\Lambda$ correspond to cellular
decompositions of a circle $S^1$. For any positive integer $n \geq
1$, we have a natural object $[n] \in \Lambda$ corresponding to the
unique decomposition with $n$ $0$-cells, called {\em vertices}, and
$n$ $1$-cells. We denote the set of vertices by $V([n])$.  Morphisms
in $\Lambda$ are homotopy classes of cellular maps satisfying
certain condition. In particular, any morphism $f:[n'] \to [n]$
induces a map $f:V([n']) \to V([n])$; moreover, it is known that for
any $v \in V([n])$, the preimage $f^{-1}(v) \subset V([n'])$ carries
a natural total order.

Now, for any monoidal category $\E$, one constructs a category
$\E^\hush$ as follows:
\begin{enumerate}
\item Objects in $\E^\hush$ are pairs $\langle [l],E_\idot \rangle$
  of an object $[l] \in \Lambda$ and a collection $E_v \in \E$, $v
  \in V([l])$ of objects in $\E$ numbered by vertices of $[l]$.
\item Morphisms from $\langle [l'],E'_\idot \rangle$ to $\langle
  [l],E_\idot \rangle$ are given by a map $f:[l'] \to [l]$ in
  $\Lambda$ and a collection of morphisms
$$
f_v:\bigotimes_{v' \in f^{-1}(v)}E'_{v'} \to E_v, \qquad v \in
V([n]),
$$
where the product is taken in the order prescribed by the total
order on $f^{-1}(v)$.
\end{enumerate}
We have a natural projection $\rho:\E^\hush \to \Lambda$ sending
$\langle [l],E_\idot \rangle$ to $[l]$. This projection is a
cofibration whose fiber over $[l] \in \Lambda$ is the product
$\E^{V([l])}$ of copies of $\E$ numbered by vertices $v \in
V([l])$. A map in $\E^\hush$ is cocartesian with respect to $\rho$
if and only if all its components $f_v$ are invertible.

To avoid size issues, assume that the monoidal category $\E$ is
small. Then trace functors from $\E$ to $\C$ form a category
$\Tr(\E,\C)$ in an obvious way, and we have the following result.

\begin{lemma}\label{hush.le}
The category $\Tr(\E,\C)$ is equivalent to the category of functors
$P^\hush:\E^\hush \to \E$ such that $P^\hush(f)$ is invertible for
any map $f$ in $\E^\hush$ cocartesian with respect to the projection
$\rho$.\endproof
\end{lemma}

The proof is in \cite[Lemma 2.5]{Ka.tr}; let us just say that the
correspondence $P^\hush \mapsto P$ simply restricts $P^\hush$ to the
fiber $\E^\hush_{[1]} \cong \E$ of the projection $\rho:\E^\hush \to
\Lambda$ over the object $[1] \in \Lambda$.

\begin{remark}\label{hush.rem}
Lemma~\ref{hush.le} has one immediate corollary. Note that the
forgetful functor $\Tr(\E,\C) \to \Fun(\E,\C)$ is faithful. Thus if
we have two functors $P^\hush_1,P^\hush_2:\E^\hush \to \C$
corresponding to trace functors $P_1,P_2 \in \Tr(\E,\C)$, and two
maps $a,a':P_1^\hush \to P_2^\hush$ such that $a=a'$ on
$\E^\hush_{[1]} \subset \E^\hush$, then Lemma~\ref{hush.le} shows
that $a=a'$ on the whole $\E^\hush$.
\end{remark}

If a monoidal category $\E$ is symmetric, for example if $\E$ is the
category of vector spaces over a field $k$, then any functor $P:\E
\to \C$ is tautologically a trace functor, with the maps
$\tau_{\idot,\idot}$ provided by the commutativity isomorphism of
the category $\E$. However, even in this case, there could be some
non-trivial trace functor structures. The basic example of such
considered in detail in \cite{Ka.tr} concern the cyclic power
functor $C_l$, $l \geq 1$ of \eqref{C.l}, and an explanation of how
it works in terms of the eqivalence of Lemma~\ref{hush.le} has been
given in \cite[Subsection 4.1]{Ka.dege}. Let us reproduce it.

Recall (e.g.\ from \cite[Appendix]{FT}) that for any integer $l \geq
1$, the category $\Lambda$ has a cousin $\Lambda_l$ described as
follows. For any $[n] \in \Lambda$, the automorphism group
$\Aut([n])$ is the cyclic group $\Z/n\Z$ generated by the clockwise
rotation $\sigma$ of the circle. If $n=ml$ is divisible by $l$, then
$\tau=\sigma^m$ generates the cyclic subgroup $\Z/l\Z \subset
\Z/n\Z$. Then objects in $[m] \in \Lambda_l$ are numbered by
positive integers $m \geq 1$, and morphisms from $[m']$ to $[m]$ are
$\tau$-equivariant morphisms from $[m'l]$ to $[ml]$ in the category
$\Lambda$. Sending $[m]$ to $[ml]$ gives a functor $i_l:\Lambda_l
\to \Lambda$. On the other hand, taking the quotient of a circle
$S^1$ by the automorphism $\tau$ and equipping it with the induced
cellular decomposition gives a functor $\pi_l:\Lambda_l \to
\Lambda$, $\pi_l([m]) = [m]$. The functor $\pi_l$ is a bifibration
with fiber $\ppt_l$. For any $m$, it also induces the natural
quotient map $\pi_l:V(i_l([m])) \to V([m])$.

Now, it was observed in \cite{Ka.dege} that for any monoidal
category $\E$ and integer $l \geq 1$, we can construct a canonical
commutative diagram
\begin{equation}\label{big.cart}
\begin{CD}
\E^\hush @<{i^{\E}_l}<< \E^\hush_l @>{\pi^{\E}_l}>> \E^\hush\\
@V{\rho}VV @VV{\rho}V @VV{\rho}V\\
\Lambda @<{i_l}<< \Lambda_l @>{\pi_l}>> \E^\hush,
\end{CD}
\end{equation}
where the square on the right-hand side is cartesian, and the
functor $i^{\E}_l$ sends $\langle [m],c_\idot \rangle$ to $\langle
i_l([m]),c^l_\idot \rangle$, with the collection $c^l_\idot$ given
by $c^l_v = c_{\pi_l(v)}$, $v \in V(i_l([m]))$. Both $\pi^\E_l$ and
$i^\E_l$ are cocartesian functors with respect to the cofibrations
$\rho$.

Assume now that $\E$ is the category of finite-dimensional
$k$-vector spaces, so that the category $\Fun(\E^\hush,k)$ is
well-defined, and consider the object
$$
C^\hush_l=\pi^{\E}_{l!}i_l^{\E *}T \in \Fun(\E^\hush,k),
$$
where $T \in \Fun(\E^\hush,k)$ is the object corresponding to the
tautological trace functor $\E \to k\amod$, $V \mapsto V$. Then for
any $V \in \E \cong \E^\hush_{[1]} \subset \E^\hush$, we have
$C^\hush_l(V) \cong C_l(V)$ by base change. Also by base change,
$C^\hush_l:\E^\hush \to k\amod$ satisfies the conditions of
Lemma~\ref{hush.le}, and corresponds to a non-trivial trace functor
structure on $C_l$.

\subsection{Constructions.}\label{tr.witt.subs}

We now observe that exactly the same construction as for the cyclic
power functor can be used to make the polynomial Witt vectors $W$
into a trace functor. Namely, let $k$ be a perfect field of
characteristic $p > 0$, fix integers $n \geq m \geq 1$, and denote
by $\E_{(n)} = W_n(k)\amod^{ff} \subset W_n(k)\amod$ the category of
finitely generated free modules over the Witt vectors ring
$W_n(k)$. We have the quotient functor $\E_{(n)} \to \E$, $V \mapsto
V/p$, where $\E = k\amod^f$ is the category of finite-dimensional
$k$-vector spaces. Since the quotient functor is monoidal, it
induces a functor $q:\E^\hush_{(n)} \to \E^\hush$. Consider the
diagram \eqref{big.cart} with $l=p^m$, and to simplify notation, let
$$
i^{(n)}_{p^m} = i^{\E_{(n)}}_{p^m}, \qquad \pi^{(n)}_{p^m} =
\pi^{\E_{(n)}}_{p^m}.
$$
Then $\pi_{p^m}^{(n)}$ is a bifibration with fiber $\ppt_{p^m}$, so
that we have the functor $\vpi_{p^m*}^{(n)}$ of \eqref{v.ga}. Denote
\begin{equation}\label{q.m}
Q_m^\hush = \vpi^{(n)}_{p^m*}i^{(n)*}_{p^m}T^{(m)} \in
\Fun(\E^\hush_{(n)},W_n(k)),
\end{equation}
where $T \in \Fun(\E^\hush_{(n)},W_n(k))$ is the tautological
functor corresponding to the embedding $\E_{(n)} \subset
W_n(k)\amod$, and $T^{(m)}$ is $T$ with the $W_n(k)$-module
structure twisted by $F^m$ as in \eqref{F.tw}. Note that by base
change, $Q_m^\hush$ satisfies the conditions of
Lemma~\ref{hush.le}. Moreover, for any $E \in \E_{(n)} =
\E^\hush_{(n)[1]} \subset \E^\hush_{(n)}$, we have $Q^\hush_m(E)
\cong \vH^0(G_m,E_{(m)})$, where $E_{(m)}$ is as in
Proposition~\ref{witt.prop}.

\begin{prop}\label{trace.prop}
\begin{enumerate}
\item For any integer $m \geq 1$, there exists an object $W^\hush_m
  \in \Fun(\E^\hush,W_m(k))$ such that $Q^\hush_m \cong
  q^*W^\hush_m$, independently of the choice of an integer $n \geq
  m$, and we have $W^\hush_m(E) \cong W_m(E)$ for any $E \in \E
  \cong \E^\hush_{[1]} \subset \E^\hush$.
\item Moreover, for any $m \geq 1$, the restriction map $R$ of
  \eqref{R.eq} extends to a map $W^\hush_{m+1} \to W^\hush_m$, and
  the Teichm\"uller representative map \eqref{T.m} extends to a map
  $T:W^\hush_1 \to W^\hush_m$.
\end{enumerate}
\end{prop}

\proof{} For \thetag{i}, note that $q$ is essentially surjective, so
that as in Proposition~\ref{witt.prop}, the issue is the morphisms:
we have to check that for two morphisms $a$, $b$ in $\E^\hush_{(n)}$
with $q(a)=q(b)$, we have $Q^\hush_m(a)=Q^\hush_m(b)$. Every
morphism $f$ in $\E^\hush_{(n)}$ decomposes as $f_v \circ f_c$,
where $f_c$ is cocartesian with respect to the projection
$\rho:\E^\hush_{(n)} \to \Lambda$, and $f_v$ is contained in a fiber
of this projection. If $q(a)=q(b)$, then in particular
$\rho(a)=\rho(b)$, so that we may assume that $a_c=b_c$, and it
suffices to check that $Q^\hush_m(a_v)=Q^\hush_m(b_m)$. Moreover,
every object $[n] \in \Lambda$ admits a morphism $f:[n] \to [1]$, so
that for any object $\wt{E} \in \E^\hush$, we have a cocartesian map
$\wt{E} \to \wt{E}_0$ with $\rho(\wt{E}_0)=[1]$. Therefore we may
further assume that $a_v$ and $b_v$ lie in the fiber $\E \cong
\E^\hush_{[1]}) \subset \E^\hush$ of the cofibration $\rho$. Then
the claim immediately follows from Proposition~\ref{witt.prop}.

For \thetag{ii}, the claim about the Teichm\"uller map $T$ is clear
from its explicit description given in the proof of
Lemma~\ref{coc.le}, so what we need to construct is the restriction
map $R$. To do this, we imitate the construction of the restriction
map $R$ of Corollary~\ref{restr.corr}. Note that the functors
$i_{p^{m+1}}$, resp.\ $\pi_{p^{m+1}}$ factor through $i_{p^m}$,
resp.\ $\pi_{p^m}$ --- we have natural functors
$$
\wt{i},\wt{\pi}:\Lambda_{p^{m+1}} \to \Lambda_{p^m}
$$
such that $i_{p^{m+1}} \cong i_{p^m} \circ \wt{i}$ and
$\pi_{p^{m+1}} \cong \pi_{p^m} \circ \wt{\pi}$. The functor
$\wt{\pi}$ is a bifibration with fiber $\ppt_p$; over an object of
$\Lambda$, it is the fibration $\ppt_{p^{m+1}} \to \ppt_{p^m}$
corresponding to the quotient map $G_{m+1} \to G_m$. Moreover, both
$\wt{\pi}$ and $\wt{i}$ lift to functors
$$
\wt{i}^{(n)},\wt{\pi}^{(n)}:\E^\hush_{(n)p^{m+1}} \to
\E^\hush_{(n)p^m}
$$
so that we have the same factorization, and $\wt{\pi}^{(n)}$ is
also a bifibration with fiber $\ppt_p$. Denote
$$
Q^{'\hush}_m =
\vpi^{(n)}_{p^{m+1}*}\wt{\pi}^{(n)*}i^{(n)*}_{p^m}T^{(m)}.
$$
Then for any $E \in \E_{(n)} = \E^\hush_{(n)[1]} \subset
\E^\hush_{(n)}$, we have $Q^{'\hush}_m(E) \cong
\vH^0(G_m,E'_{(m)})$, where $E'_{(m)}$ is as in
Proposition~\ref{witt.prop}, so that by the same argument as in
\thetag{i}, we have
$$
Q^{'\hush}_m = q^*W^{'\hush}_m
$$
for some $W^{'\hush}_m \in \Fun(\E^\hush,W_n(k))$. Moreover, as in
\eqref{r.eq}, we have a natural map
$$
r:W^{'\hush}_m \to W^\hush_m,
$$
so that it suffices to construct an isomorphism $W^{'\hush}_m \cong
W^\hush_{m+1}$. By Lemma~\ref{hush.le}, this is equivalent to
proving that the isomorphism \eqref{c.eq} of
Corollary~\ref{restr.corr} is compatible with the trace functor
structures --- that is, commutes with the maps $\tau_{\idot,\idot}$.

However, the category $\Gamma$ of finite sets is also a monoidal
category, with the monoidal structure given by cartesian product,
and we have a natural mono\-idal functor $\nu:\Gamma \to \E$ sending
a finite set $S$ to the vector space $k[S]$. Then
$\nu^*W^\hush_{m+1}$, $\nu^*W^{'\hush}_m$ correspond to trace
functors from $\Gamma$ to $W(k)\amod$, and since $\nu$ is
essentially surjective, it suffices to check that \eqref{c.eq}
commutes with $\tau_{\idot,\idot}$ after restricting to $\Gamma$. In
other words, it suffices to prove that it extends to an isomorphism
$\nu^*W^{'\hush}_m \cong \nu^*W^\hush_{m+1}$. Moreover, the pullback
functor $q^*$ is fully faithful, and $\nu$ factors through $q$ by
means of a monoidal functor $\wt{\nu}:\Gamma \to
\E_{(n)}$. Therefore it suffices to extend \eqref{c.eq} to an
isomorphism $\eps:\wt{\nu}^*Q^{'\hush}_m \to
\wt{\nu}^*Q^\hush_{m+1}$. But by base change,
$\wt{\nu}^*Q^{'\hush}_m$ and $\wt{\nu}^*Q^\hush_{m+1}$ are given by
$$
\wt{\nu}^*Q^{'\hush}_m = \vpi_{p^{m+1}!}^\Gamma \wt{\pi}^{\Gamma
  *}\wt{T}^{(m)},
\qquad 
\wt{\nu}^*Q^\hush_{m+1} = \vpi_{p^{m+1}!}^\Gamma
\wt{i}^{\Gamma *}\wt{T}^{(m)},
$$
where we denote $\wt{T}^{(m)} = i^{\Gamma
  *}_{p^n}\wt{\nu}^*T^{(m)}$. The canonical maps $c_S$, $S \in
\Gamma$ of \eqref{c.S} together give a map
\begin{equation}\label{adm.eq}
\wt{\pi}^{\Gamma *}\wt{T}^{(m)} \to \wt{i}^{\Gamma *}\wt{T}^{(m)}.
\end{equation}
This map induces a map $\eps:\wt{\nu}^*Q^{'\hush}_m \to
\wt{\nu}^*Q^\hush_{m+1}$ whose restriction to $\Gamma =
\Gamma^\hush_{[1]} \subset \Gamma^\hush$ is exactly the isomorphism
\eqref{c.eq}. Since both $Q^{'\hush}_m$ and $Q^\hush_{m+1}$ send
cocartesian maps in $\E^\hush$ to invertible maps, and every object
in $\Gamma^\hush$ admits a cocartesian map to an object in
$\Gamma^\hush_{[1]}$, the map $\eps$ must be an isomorphism
everywhere.
\endproof

As a corollary of Proposition~\ref{trace.prop}, we see that the
$m$-truncated polynomial Witt vectors functors $W_m$, $m \geq 1$
have natural trace functors structures at least if we restrict them
to finite-dimensional $k$-vector spaces, and the restriction maps
$R:W_{m+1} \to W_m$ together with the Teichm\"uller maps $T:W_1 \to
W_m$ are trace functor maps. Since there functors commute with
filtered colimits, both statements immediately extend to all
$k$-vector spacs. The inverse limit $W =
\lim_{\overset{R}{\gets}}W_m$ then also has a structure of a trace
functor, and the Teichm\"uller map \eqref{T} is a trace functor map.

\subsection{Extensions.}

Let us now prove that the trace functor structure on $W$ is
compatible with the two other structures it has --- namely, the
structure of a $G$-Mackey functor, and the pseudotensor structure of
Proposition~\ref{mult.prop}. As it happens, the proofs are quite
straightforward, and the main issue is formulating the exact meaning
of compatibility.

For the pseudotensor structure, this is also straightforward. For
any unital monoidal category $\E$, the category $\E^2 = \E \times
\E$ is also unital monoidal with respect to the coordinatewise
monoidal structure. For any trace functor $P$ from $\E$ to a
category $\C$, we have a natural trace functor $P^2:\E^2 \to
\C^2$. If $\E$ is symmetric, then the product functor $m_{\E}:\E^2
\to \E$ is a monoidal functor, so that it defines a natural functor
$$
m^\hush_{\E}:\E^{2\hush} \to \E^\hush
$$
cocartesian over $\Lambda$. Moreover, the unit object $1 \in \E$
extends to a natural cocartesian section $1_\Lambda:\Lambda \to
\E^\hush$ of the projection $\rho^\E:\E^\hush \to \Lambda$ such that
$1_\Lambda([n]) = 1^{V([n])} \in \E^{V([n])}$, $[n] \in \Lambda$.

\begin{defn}\label{tr.mult.def}
A {\em pseudotensor structure} on a trace functor $P$ from a
symmetric monoidal category $\E$ to a monoidal category $\C$ with
the product functor $m_{\C}:\C^2 \to \C$ is given by functorial maps
\begin{equation}\label{eps.mu.la}
\eps:1 \to P^\hush(1_\Lambda), \qquad \mu:m_{\C} \circ P^{2\hush} \to
P^\hush \circ m^\hush_\E
\end{equation}
that are associative, commutative and unital in the obvious sense.
\end{defn}

As for the $G$-Mackey structure, one immediate observation is that
finite $G$-orbits and objects $[n] \in \Lambda$ have one thing in
common: their automorphisms form finite cyclic groups. We emphasize
this similarity by using the notation $[p^m]$, $m \geq 0$ for finite
$G$-orbits in Subsection~\ref{cycl.mack.subs}. In fact, we have a
natural functor $\delta:I \to \Lambda$, $\delta([p^m]) = [p^m]$,
where as in Subsection~\ref{cycl.mack.subs}, $I$ is the groupoid of
finite $G$-orbits and their isomorphisms. Moreover, this functor
fits into a commutative diagram
\begin{equation}\label{del.dia}
\begin{CD}
I @<{\pi}<< I_p @>{i}>> I\\
@V{\delta}VV @VVV @VV{\delta}V\\
\Lambda @<{\pi_p}<< \Lambda_p @>{i_p}>> \Lambda,
\end{CD}
\end{equation}
where the squares are cartesian, and $I_p$, $i$, $\pi$ are as in
Subsection~\ref{cycl.mack.subs}. One is tempted to think that one
can add morphisms of $G$-orbits to morphisms already in $\Lambda$ to
obtain a bigger category $\LR \supset \Lambda$, and then try to
imitate the definition of Mackey functors to obtain some sort of a
category $Q\LR$. This is indeed possible; the resulting category is
the {\em cyclotomic category} that appeared e.g.\ in
\cite{ka.cyclo}, and when the theory is fully developed, it provides
a notion of a ``cyclotomic Mackey functor''. However, doing all this
properly requires a lot of work and a lot of space, and we will
return to it elsewhere. For the purposes of the present paper, we
simply use an explicit description of $G$-Mackey functors given in
Lemma~\ref{cycl.mack.le}, and introduce the following somewhat {\em
  ad hoc} definition.

\begin{defn}\label{tr.FV.def}
An {\em $FV$-structure} on a trace functor $P$ from a monoidal
category $\E$ to the category $A\amod$ of modules over a ring $A$ is
a pair of maps
$$
\begin{CD}
P^\hush \circ i^\E_p @>{V}>> P^\hush \circ \pi^\E_p @>{F}>> P^\hush
\circ i^\E_p
\end{CD}
$$
such that $F \circ V:P^\hush \circ i^E_p \to P^\hush \circ i^\E_p$
is equal to the natural map $\tr^\dg_{\pi^\E_p}$ of \eqref{wt.tr}.
\end{defn}

Here $i^\E_p,\pi^\E_p:\E^\hush_p \to E^\hush$ are the functors of
\eqref{big.cart}, and we recall that the trace map
$\tr^\dg_{\pi^E_p}$ is well-defined even if the category $\E$ is
not small. Note that for any object $E \in \E = \E^\hush_{[1]}
\subset \E^\hush$, we have a natural functor $E^\hush:I \to
\E^\hush$ given by
\begin{equation}\label{wt.E.hush}
E^\hush([p^m]) = i^{\E}_{p^m}(\wt{E}),
\end{equation}
and we have $\rho^\E \circ E^\hush \cong \delta$. Therefore by
\eqref{del.dia}, for any trace functor $P^\hush:\E \to A\amod$
equipped with an $FV$-structure, and for any object $E \in \E$, the
composition $P^\hush \circ E^\hush$ with the induced maps $V$, $F$
satisfies the conditions of Lemma~\ref{cycl.mack.le}~\thetag{i} and
defines an $A$-valued $G$-Mackey functor.

\smallskip

Finally, assume that the monoidal category $\E$ is symmetric, the
ring $A$ is commutative, and we are given a trace functor $P$ from
$\E$ to $A\amod$ that has both a pseudotensor structure $\langle
\eps,\mu \rangle$ in the sense of Definition~\ref{tr.mult.def} and
an $FV$-structure $\langle V,F \rangle$ in the sense of
Definition~\ref{tr.FV.def}.

\begin{defn}\label{comp.def}
The structures $\langle \eps,\mu \rangle$ and $\langle V,F \rangle$
on a trace functor $P$ from $\E$ to $A\amod$ are {\em compatible} if
we have
\begin{gather}\label{FV.m}
\mu \circ (F \times F) = F \circ \mu,\\
\mu \circ (V \times \id) = V \circ \mu \circ (\id \times F),\\
\mu \circ (\id \times V) = V \circ \mu \circ (F \times \id).
\end{gather}
\end{defn}

We can now formulate our result about Witt vectors. Note that a
pseudotensor structure on a trace functor $P$ from $\E$ to $\C$
gives a pseudotensor structure on the underlying functor $P:\E \to
\C$ by restriction to the fiber $\E = \E^\hush_{[1]} \subset
\E$. Let $\E=k\amod$, the category of vector spaces over our perfect
ring $k$ of characteristic $p$, and let $\C=W(k)\amod$, the category
of modules over the Witt vectors ring $W(k)$.

\begin{prop}
The Witt vectors trace functor $W$ from $\E = k\amod$ to abelian
groups of Subsection~\ref{tr.witt.subs} has a natural $W(k)$-linear
pseudotensor structure $\langle \eps,\mu \rangle$ and an
$FV$-structure $\langle V,F \rangle$ such that
\begin{enumerate}
\item the two structures are compatible in the sense of
  Definition~\ref{comp.def},
\item $\langle \eps,\mu \rangle$ restricts to the pseudotensor
  structure of Proposition~\ref{mult.prop} on the fiber $k\amod =
  k\amod^\hush_{[1]} \subset k\amod^\hush$, and
\item for any object $E \in \E$, the functor $W^\hush \circ
  E^\hush:I \to W(k)\amod$ with the maps $V$, $F$ corresponds to the
  extended Witt vectors Mackey functor $\wt{W}(E)$ of
  \eqref{ext.W.eq} under the equivalence of
  Lemma~\ref{cycl.mack.le}.
\end{enumerate}
Moreover, for any integer $n \geq 1$, the maps $V$, $F$ of the
$FV$-structure induce $W_{n+1}(k)$-linear functorial maps
\begin{equation}\label{FV.n}
\begin{CD}
W^{(1)\hush}_n \circ i_p^{k\amod} @>{V}>> W_{n+1}^\hush \circ \pi_p^{k\amod}
@>{F}>> W_n^{(1)\hush} \circ i_p^{k\amod},
\end{CD}
\end{equation}
and the pseudotensor structure $\langle \eps,\mu \rangle$ induces a
pseudotensor structure on the quotient $W_n$ of the trace functor
$W$.
\end{prop}

\proof{} As in the proofs of Proposition~\ref{mult.prop} and
Proposition~\ref{trace.prop}, it suffices to prove everything after
replacing $\E$ with the category $k\amod^f$ of finite-dimensional
$k$-vector spaces, and then passing to the categories $\E_{(n)} =
W_n(k)\amod^{ff}$, $n \geq 1$ by means of the fully faithful
pullback functor $q^*$. Thus to construct a pseudotensor structure
on $W$, it suffices to construct a system of pseudotensor structures
on trace functors $Q_n=q^*W_n$, $n \geq 1$ compatible with the
restriction maps $R$. Fix such an integer $n$, and to simplify
notation, let $m^\hush = m^\hush_{\E_{(n)}}$,
$i_{p^n}=i_{p^n}^{\E_{(n)}}$, $\pi_{p^n}=\pi_{p^n}^{\E_{(n)}}$. Then
since the tautological embedding $\E_{(n)} \subset W_n(k)\amod$ is a
tensor functor, we have
$$
m^{\hush *}T^{(n)} \cong T^{(n)} \boxtimes T^{(n)},
$$
and since $i_{p^n} \circ m^\hush \cong m^\hush \circ i_{p^n}$, this
gives an isomorphism
$$
m^{\hush *}i_{p^n}^*T^{(n)} \cong i_{p^n}^*T^{(n)} \boxtimes
i_{p^n}^*T^{(n)}.
$$
Therefore to obtain a pseudotensor structure on $Q^\hush_n =
\vpi_{p^n*}i_{p^n}^*T^{(n)}$, it suffices to recall that the functor
$\vpi_{p^n*}$ is pseudotensor by Lemma~\ref{vpi.le}. On the fiber
$\E = E^\hush_{[1]} \subset \E^\hush$, $Q^\hush_n$ restricts to the
functor $Q_m$, so that the compatibility statement \thetag{ii}
immediately follows from the characterisation of the pseudotensor
structure on $W_n$ given in Subsection~\ref{mult.subs}. Then
Proposition~\ref{mult.prop} implies that our pseudotensor structures
agree with the restriction maps on $\E^\hush_{[1]} \subset
\E^\hush$, and as explained in Remark~\ref{hush.rem},
Lemma~\ref{hush.le} shows that they agree everywhere.

To obtain the $FV$-structure, note that for every $n \geq 1$, we
have
$$
i_p^*Q^\hush_n \cong \vpi_{p^n*}i_{p^{n+1}}^*T^{(n)}
$$
by base change. Simplify notation further by writing $\wt{T} =
i_{p^{n+1}}^*T^{(n+1)}$, $\pi = \pi_{p^{n+1}}$, and let $\pi'=\pi_p$,
$\pi''=\pi_{p^n}$. Then \eqref{lr.pi} provides natural maps
$l^\dg_\pi$, $r^\dg_\pi$, and evaluating them at $\wt{T}$, we obtain
morphisms
$$
V = l^\dg_\pi(\wt{T}):i_p^*Q^{(1)\hush}_n \to \pi_p^*Q^{\hush}_{n+1},
\qquad F = r^\dg_\pi(\wt{T}):\pi_p^*Q^\hush_{n+1} \to
i_p^*Q^{(1)\hush}_n
$$
whose composition $F \circ V$ is exactly as required in
Definition~\ref{tr.FV.def} by virtue of \eqref{fact}. Moreover,
since $q^*$ is fully faithful, the maps $F$ and $V$ descend to Witt
vectors trace functors $W^\hush_\idot$ and give the diagram
\eqref{FV.n}.

Now observe that if choose an object $E \in \E$ and restrict
$W_n^\hush$ and $W_{n+1}^\hush$ to the category $I$ via the functor
\eqref{wt.E.hush}, then by \eqref{FV.lr}, our maps $V$ and $F$
restrict exactly to the maps $V$ and $F$ corresponding to the Mackey
functor $\wt{W}_{n+1}(E)$ under the equivalence of
Lemma~\ref{cycl.mack.le}. This implies in particular that the maps
$V$ and $F$ for different integers $n$ agree with the restriction
maps $R$ on $\E = \E^\hush_{[1]} \subset \E^\hush$, and then by
Lemma~\ref{hush.le}, they must agree everywhere. Therefore we can
pass to the inverse limit and obtain an $FV$-structure on the trace
functor $W^\hush$. Moreover, this $FV$-structure satisfies the
compatibility condition \thetag{iii}.

To finish the proof, it remains to notice that the remaining
condition \thetag{i} amounts to checking \eqref{FV.m}, and this can
be done pointwise, that is, after evaluation at an arbitrary object
$\wt{E} \in \E^\hush$. Moreover, by Lemma~\ref{hush.le}, it suffices
to consider objects $\wt{E} \in \E^\hush_{[1]} \subset \E^\hush$,
and by \thetag{iii}, the statement then immediately follows from
Lemma~\ref{cycl.mack.le}~\thetag{ii}.
\endproof

{\small\noindent
Affiliations (in the precise form required for legal reasons):
\begin{enumerate}
\renewcommand{\labelenumi}{\arabic{enumi}.}
\item Steklov Mathematics Institute, Algebraic Geometry Section
  (main affiliation).
\item Laboratory of Algebraic Geometry, National Research University
Higher\\ School of Economics.
\item Center for Geometry and Physics, Institute for Basic
  Science (IBS), Pohang, Korea.
\end{enumerate}}

\noindent
{\em E-mail address\/}: {\tt kaledin@mi.ras.ru}

\end{document}